\documentclass[12pt,reqno]{amsart}
\newtheorem{thm}{Theorem}[section]
\newtheorem{defi}{Definition}[section]

\newtheorem{cor}{Corollary}[section]
\newtheorem{pr}{Proposition}[section]
\theoremstyle{definition}
\newtheorem{rem}{Remark}[section]
\newtheorem{note}{Note}[section]

\newcommand{\be}{\begin{equation}}
\newcommand{\ee}{\end{equation}}
\newcommand{\bea}{\begin{eqnarray}}
\newcommand{\eea}{\end{eqnarray}}
\newcommand{\beb}{\begin{eqnarray*}}
\newcommand{\eeb}{\end{eqnarray*}}
\usepackage{amssymb,amsfonts,amsthm,setspace,indentfirst,mathrsfs}
\usepackage{minibox}
\usepackage{enumitem}

\numberwithin{equation}{section}
\begin{document}
%
\title[Vaidya-Bonner-de Sitter spacetime]{Symmetry and pseudosymmetry properties of Vaidya-Bonner-de Sitter spacetime}
\author[A. A. Shaikh, S. K. Hui, M. Sarkar \& V. A. Babu]{Absos Ali Shaikh$^{* 1}$, Shyamal kumar Hui$^2$, Mousumi Sarkar$^3$ \& V. Amarendra Babu $^4$}
\date{\today}
\address{\noindent\newline$^{1,2,3}$ Department of Mathematics,
	\newline University of Burdwan, 
	\newline Golapbag, Burdwan-713104,
	\newline West Bengal, India} 

\address{\noindent\newline$^{4}$ Department of Mathematics,
\newline Achraya Nagarjuna University,
\newline Nagarjuna Nagar, Guntur-522510,
\newline Andhra Pradesh, India}
\email{aask2003@yahoo.co.in$^1$, aashaikh@math.buruniv.ac.in$^1$}
\email{skhui@math.buruniv.ac.in$^2$}
\email{mousrkr13@gmail.com$^3$}
\email{amarendrab4@gmail.com$^4$}

\dedicatory{}

\begin{abstract}
The primary focus of the current study is to explore the geometrical properties of the Vaidya-Bonner-de Sitter (briefly, VBdS) spacetime, which is a generalization of Vaidya-Bonner spacetime, Vaidya spacetime and Schwarzschild spacetime. In this study we have shown that the VBdS spacetime describes various types of pseudosymmetric structures, including pseudosymmetry due to conformal curvature,  conharmonic curvature and other curvatures. Additionally, it is shown that such a spacetime is $2$-quasi-Einstein, Einstein manifold of level $3$, generalized Roter type, and that conformal $2$-forms are recurrent. The geometric features of the Vaidya-Bonner spacetime, Vaidya spacetime, and Schwarzschild spacetime are obtained as a particular instance of the main determination. It is further established that the VBdS spacetime admits almost Ricci soliton and almost $\eta$-Yamabe soliton with respect to non-Killing vector fields. Also, it is proved that such a spacetime possesses generalized conharmonic curvature inheritance. It is interesting to note that in the VBdS spacetime the tensors $Q(T,R)$, $Q(S,R)$ and $Q(g,R)$ are linearly dependent. Finally, this spacetime is compared with the Vaidya-Bonner spacetime with respect to their admitting geometric structures,  viz., various kinds of symmetry and pseudosymmetry properties. 
 
\end{abstract}

%
\maketitle
%

\section{\bf Introduction}\label{intro}

Let $(\mathscr{V}_n,g)$, $n=dim (\mathscr{V})\geq 3$, be a smooth and connected manifold, endowed with the semi-Riemannian metric $g$ of signature $(\chi, n-\chi)$, $0\leq \chi \leq n$ (throughout the paper we use this notation).  The signature is called Lorentzian if $\chi=1$ or $\chi=n-1$,  and $\mathscr{V}_n$ is a spacetime if it admits Lorentzian signature with $n=4$. Let $R$, $S$, $\bar{\kappa}$, $C$, $P$, $har(R)$ and $cir(R)$ be respectively the Riemann-Christoffel,  Ricci,  scalar, Weyl conformal,  projective,  conharmonic and concircular curvature tensor of $(\mathscr{V}_n,g)$ along with the Levi-Civita connection  $\nabla$. 

\indent  Curvature elegantly discloses the geometric properties of a manifold, which plays a crucial role in grasping the shape of that manifold. The symmetry of a manifold is determined locally by restricting the curvature tensor $R$, as, in $1926$, Cartan \cite{Cart26} established the idea of the local symmetry of a manifold defined by $\nabla R=0$. 
Again, the concept of semisymmetry, which is a generalization of local symmetry, was also introduced by Cartan \cite{Cart46}. Later, Szab\'o \cite{Szab82, Szab84, Szab85} classified the notion of semisymmetry for the Riemannian manifolds.
After that, generalizing such notions in different directions, Cartan's works were further developed by imposing more weakened curvature restrictions.  Such generalized notions of symmetries are as follows: the recurrence of Ruse \cite{Ruse46, Ruse49a, Ruse49b} (also, see  \cite{Walk50}), the super, hyper and weakly generalized recurrence by Shaikh et al. \cite{SAR13, SKA18, SP10, SR10, SR11, SRK18, SRK17},  curvature $2$-forms of recurrent manifolds of Besse \cite{Bess87} (also, see \cite{LR89, MS12a, MS13a, MS14}), weak symmetry of Tam\'{a}ssy and Binh \cite{TB89, TB93},  pseudosymmetric manifolds of Adam\'{o}w and Deszcz \cite{AD83} and pseudosymmetric manifolds of Chaki \cite{Chak87, Chak88}. Recently,  many authors (see e.g. \cite{ADEHM14, DDHKS00, DK99, HV07, HV09, Kowa06, SAA18, SAAC20, SAACD_LTB_2022, SAD_pgm_2023,  SDKC19,   SHS_warped, SK16srs}) investigated different spacetimes to examine whether such spacetimes admit the structures of symmetry, weak symmetry, recurrency and pseudosymmetry. Additionally, Mike\v s et al. \cite{MIKES76, MIKES88, MIKES92, MIKES96, MS94, MSV15, MVH09} have thoroughly investigated the idea of geodesic mappings of many types of symmetric Riemannian manifolds.    We note that the notions of pseudosymmetry introduced independently by Deszcz \cite{Desz92} and  Chaki \cite{Chak87} are completely different. \\


In $1982$, Hamilton \cite{Hamilton1982} proposed the concept of Ricci flow while studying $3$-dimensional compact manifolds with positive Ricci curvature. The Ricci flow is a process of   changing the Riemannian metric over time.  Ricci solitons, which are natural generalizations of Einstein metrics \cite{Bess87, Brink1925, S09, SYH09}, are the self-similar solutions of the Ricci flow. The term ``Ricci soliton" has been generalized in several ways, including almost Ricci soliton,  almost $\eta$-Ricci soliton,  and $\eta$-Ricci soliton. 
If a semi-Riemannian manifold $\mathscr{V}_n$ possesses the relation
$$\frac{1}{2}\pounds_\xi g+S-\delta g=0,$$
 where $\delta$ is the non-constant smooth function, then $\mathscr{V}_n$ is called an almost Ricci soliton (\cite{Pigola2011}). Again, if $\delta$= constant, then it is known as a Ricci soliton, where $\pounds_\xi$ is the Lie derivative concerning the soliton vector field $\xi$. We note that it is steady,  expanding, or shrinking if $\delta=0$, $\delta<0$, or $\delta>0$ respectively.  
Again, if  a non-zero $1$-form $\eta$ on $\mathscr{V}_n$ possesses the relation 
$$\frac{1}{2}\pounds_\xi g+S-\delta g+\nu (\eta\otimes\eta)=0,$$  where $\delta,\nu$ are non-constant smooth functions, then $\mathscr{V}_n$ is called an almost $\eta$-Ricci soliton \cite{Blaga2016}. Again, $\mathscr{V}_n$ is called $\eta$-Ricci soliton \cite{Cho2009}  if $\delta, \nu$ are constants.\\  

Another geometric flow, known as the Yamabe flow, was  introduced by Hamilton \cite{Hamilton1988} simultaneously with the Ricci soliton. G\"uler and Cr\'{a}\c{s}mare\v{a}nu \cite{Guler2019} recently developed another geometric flow known as Ricci-Yamabe flow as a linear combination of both the Ricci and Yamabe flow. The self-similar solutions of a Ricci-Yamabe (or Yamabe) flow are known as Ricci-Yamabe (or Yamabe) solitions. Again, if $\mathscr{V}_n$ realizes the condition 
$$\frac{1}{2}\pounds_\xi g+\beta_1 S+\left(\delta-\frac{1}{2}\beta_2\bar{\kappa}\right) g=0,$$
 then it is called a Ricci-Yamabe soliton  \cite{Siddiqi2020}, where $\beta_1$, $\beta_2$, $\delta$ are constants. Also, it is called an almost Ricci-Yamabe soliton \cite{Siddiqi2020} if  $\beta_1$, $\beta_2$, $\delta$ are  non-constant smooth functions on $\mathscr{V}_n$. We note that if $\beta_1=0$, $\beta_2=2$, then it becomes a Yamabe soliton.  Again, if a non-zero $1$-form $\nu$ satisfies
$$\frac{1}{2}\pounds_\xi g+\beta_1 S+\left(\delta-\frac{1}{2}\beta_2\bar{\kappa}\right) g+\nu \eta\otimes\eta=0,$$
  then it is known as an almost $\eta$-Ricci-Yamabe soliton \cite{Siddiqi2020}  with the non-constant smooth functions $\beta_1$, $\beta_2$, $\delta$, $\nu$. In addition, if $\beta_1$, $\beta_2$, $\delta$, $\nu$ are constants, then it is known as a $\eta$-Ricci-Yamabe soliton \cite{Siddiqi2020}.
Over the past three decades, a large number of research articles on Ricci soliton, Yamabe soliton, and their generalizations (see, \cite{Ahsan2018, AliAhsan2013, AliAhsan2015, SDAA_LCS_2021} and the references therein) have been published; as a result, it is quite an active area of research in the realm of differential geometry.

Again, imposing symmetry is a critical tool to create gravitational potentials that fulfill Einstein field equations, indicating that geometrical symmetries are indispensable for the investigation of the theory of general relativity. If the Lie derivative of a certain tensor field disappears with respect to a non-zero vector field and the vanishing Lie derivative exhibits geometric symmetry, then a geometric quantity is preserved along that vector field. Such symmetries can be described as motion, curvature collineation, Ricci collineation, etc. The significance of curvature collineation in general relativity was thoroughly investigated by Katzin {et al.} \cite{KLD1969, KLD1970}. In $1992$, Duggal [1] established the concept of curvature inheritance, which expanded the idea of curvature collineation for the curvature tensor of type (1,3). In the last three decades, a large number of publications (see, \cite{ Ahsan1977_231, Ahsan1977_1055, Ahsan1978, Ahsan1987, Ahsan1995, Ahsan1996, Ahasan2005, AA2012, AhsanAli2014,  AH1980, AliAhsan2012, ShaikhDatta2022})  about the studies of these types of symmetries have appeared in the literature. Shaikh and Datta [2] recently developed the idea of generalized curvature inheritance, which is an extension of curvature collineation and curvature inheritance for the type of (0,4) curvature tensor.

    The VBdS spacetime \cite{VBdS_1987} is an exact solution of the Einstein field equations with a non-zero cosmological constant, which represents the gravitational field around a spherically symmetric mass distribution.
The line element of the VBdS spacetime  in advanced Eddington-Finkelstein time coordinates $(t,r,\theta,\phi)$ is given as
    \begin{eqnarray}\label{VBdS}
	ds^2=\left(1-\frac{2m(t)}{r}+\frac{q^2(t)}{r^2}-\frac{\lambda r^2}{3}\right)dt^2-2dtdr-r^2(d\theta^2+\sin^2\theta d\phi^2),
    \end{eqnarray}
     where $\lambda$ is the cosmological constant, $t$ is the advanced time coordinate, $r$ is the radial coordinate, and both the mass $m(t)$ and charge $q(t)$ of the body depend on time. The VBdS spacetime is a generalization of Vaidya-Bonner spacetime, Vaidya spacetime, and Schwarzschild spacetime. In fact, 
    \begin{enumerate}[label=(\alph*)]
        \item  if $\lambda=0$ in (\ref{VBdS}), then it becomes Vaidya-Bonner spacetime,


        \item if $\lambda=0=q(t)$ and the mass of the body is time independent in (\ref{VBdS}), it turns into Vaidya spacetime.

        \item If $m(t)=$ constant and $\lambda=0=q(t)$, the VBdS spacetime (\ref{VBdS}) reduces to the Schwarzschild spacetime.
        \end{enumerate}
During the last two decades, several physical phenomenons, such as quantum effect \cite{VBdS_1999}, thermal radiation of scalar particles \cite{VBdS_2001}, entropy \cite{VBdS_2012}, quantum non-thermal effect \cite{VBdS_2013}, Lorentz-breaking theory \& tunneling radiation correction \cite{VBdS_2022}, Hawking radiation \cite{VBdS_2007_Hawking}  have been investigated on the VBdS spacetime by various physicists.

\indent Although numerous researchers have found the curvature-restricted geometric structures of many spacetimes, the geometric properties of VBdS spacetime have not been fully investigated.  The primary objective of the work is to investigate the curvature aspects of the VBdS spacetime. The VBdS spacetime is shown to possesses several types of pseudosymmetric properties, including pseudosymmetry due to conformal curvature, pseudosymmetry due to conharmonic curvature, and the difference tensor $R \cdot R - Q(g, C)$ is linearly dependent on $Q(S, R)$. Even, the VBdS spacetime is neither semisymmetric nor pseudosymmetric. Also, the nature of such spacetime is of the generalized Roter type, Einstein manifold of level $3$ and $2$-quasi Einstein manifold.  The geometric features of the Vaidya-Bonner spacetime, Vaidya spacetime, and Schwarzschild spacetime are deduced as a particular instance of the main result.  Here, we have also found that this spacetime admits almost Ricci soliton and almost $\eta$-Yamabe soliton with respect to non-Killing vector fields. Also, it is proved that such a spacetime possesses generalized conharmonic curvature inheritance.
Finally, the geometric structures of this spacetime and the Vaidya-Bonner spacetime are compared.

\indent The article is organized as follows: Section $2$ discusses various types of curvatures as well as several significant curvature criteria that define the geometric structures of a manifold. Section $3$ is devoted to the main investigation on VBdS spacetime. The nature of the Ricci and Ricci-Yamabe solitons admitted by the VBdS spacetime is explained in Section $4$. In Section $5$, we have studied the geometric nature of energy momentum tensor of the VBdS spacetime. Comparing the admitting geometric structures of the VBdS spacetime with the Vaidya-Bonner spacetime, various kinds of symmetry and the Ricci soliton have discussed in Section $6$.  \\



\section{\bf Preliminaries}

This part aims to discuss various geometric structures that occur with restrictions on curvatures and their covariant derivatives, which are useful to explain the symmetry of the VBdS spacetime and have related geometric meanings. Additionally, this section provides explanations of the concepts of motion, curvature inheritance, Ricci inheritance, curvature collineation and Ricci collineation, etc.

Two symmetric tensors $X$ and $Z$ of type $(0, 2)$ combine to form the Kulkarni-Nomizu product $X \wedge Z$, which is defined as 
(see, \cite{DGHS11,DHJKS14, Glog02, G08, Kowa06, SC21, SK19}):
$$(X\wedge Z)_{efst}=2X_{e[t}Z_{s]f} +2X_{f[s}Z_{t]e},$$
where the antisymmetrization of index pairs is denoted by $[.]$. Now, the Riemann curvature of type $(1,3)$, conharmonic curvature of type $(1,3)$, concircular curvature of type $(1,3)$, Weyl conformal curvature of type $(1,3)$, and projective curvature of type $(1,3)$ are given as follows:
\bea
R^e_{fst}&=& 2\left(\Gamma^m_{f[s}\Gamma^e_{t]m} + \partial_{[t}\Gamma^e_{s]f}\right),\nonumber\\
har(R)^e_{fst}&=& R^e_{fst} - \frac{2}{n-2}\left( J^e_{[f}g_{s]t} + \delta^e_{[f}S_{s]t}\right),\nonumber \\
cir(R)^e_{fst}&=& R^e_{fst} - \bar\kappa\frac{2}{n(n-1)}\delta^e_{[f}g_{s]t}, \nonumber \\
C^e_{fst}&=& R^e_{fst}+\frac{2}{n-2}\left( \delta^e_{[f}S_{s]t} +J^e_{[f}g_{s]t}\right) -\bar\kappa\frac{2}{(n-1)(n-2)}\delta^e_{[f}g_{s]t} ,\nonumber \\
P^e_{fst}&=& R^e_{fst} -\frac{2}{n-1}\delta^e_{[f}S_{s]t}, \nonumber
\eea
where $\partial_e=\frac{\partial}{\partial x^e}$,  $J^f_s$ is the $(1,1)$-type Ricci curvature,  and  $\Gamma^f_{st}$ are the connection coefficients.  By simply lowering the indices with the help of the metric $g_{ef}$, one may obtain the $(0,4)$ rank tensors $R_{efst}$, $har(R)_{efst}$, $cir(R)_{efst}$, $C_{efst}$ and $P_{efst}$ which are given as follows:
\bea
R_{efst}&=& g_{e\alpha}(\partial_t \Gamma^\alpha_{fs}-\partial_s \Gamma^\alpha_{fs}+\Gamma^\beta_{fs}\Gamma^\alpha_{\beta t}-\Gamma^\beta_{ft}\Gamma^\alpha_{\beta s}) , \nonumber \\
har(R)_{efst}&=& R_{efst} - \frac{1}{n-2} (g\wedge S)_{efst} , \nonumber \\
cir(R)_{efst}&=& R_{efst} - \frac{\bar \kappa}{2n(n-1)} (g\wedge g)_{efst}, \nonumber \\
C_{efst}&=& R_{efst}-\frac{1}{n-2}(g\wedge S)_{efst}+\frac{\bar\kappa}{2(n-1)(n-2)}(g\wedge g)_{efst} ,\nonumber \\
P_{efst}&=& R_{efst} -\frac{1}{n-1}(g_{et}S_{fs}-g_{ft}S_{es}). \nonumber
\eea

Let $W$ be a $(0,m)$, $m\geq 1$, tensor. We also define the $(0,m+2)$ rank tensors $L \cdot W$ (see, \cite{DG02, DGHS98, DH03, SDHJK15, SK14}) and $Q(\beta, W)$ (see, \cite{DGPSS11,SDHJK15,Tach74}) as follows:
\beb
(L\cdot W)_{b_1b_2\cdots b_mrs}&=&-\left[ L^\alpha_{rsb_1}W_{\alpha b_2\cdots b_m}+ \cdots + L^\alpha_{rsb_m}W_{b_1\cdots \alpha}\right], \\
Q(\beta,W)_{b_1b_2\cdots b_m rs}&=&\beta_{sb_1}W_{rb_2\cdots b_m}+ \cdots + \beta_{sb_m}W_{b_1b_2\cdots r}\\ 
&-& \beta_{rb_1}W_{sb_2\cdots b_m}- \cdots - \beta_{rb_m}\beta_{b_1b_2\cdots s},
\eeb

where $\beta_{ef}$ is a symmetric $(0,2)$ rank tensor and $L_e^{fst}$ is a $(1,3)$ rank tensor.

\begin{defi} \cite{AD83, Cart46,  Desz92, Desz93, DGHZ15, DGHZ16,SAAC20N, SK14, SKppsnw, Szab82, Szab84, Szab85} 
Let $\mathscr{V}_n$ be a semi-Riemannian manifold. If in $\mathscr{V}_n$, the relation $L\cdot W=0$ holds, then it is termed as a $W$-semisymmetric type manifold due to $L$.  Also, $\mathscr{V}_n$ is said to be $W$-pseudosymmetric type manifold due to $L$ if the relation $L\cdot W=\mathscr{F}_L Q(\beta, W)$ holds, for a smooth function $\mathscr{F}_L$ on the set $\{ x\in M: Q(\beta,W)\neq 0 \ at \ x \}$ (i.e., the tensors $L\cdot W$ and $Q(\beta,W)$ are linearly dependent).
\end{defi}
If we substitute $L$ = $R$ and $W$ = $R$ (resp., $P$, $har(R)$, $cir(R)$, $C$ and $S$) in the above definition, then the $W$-semisymmetric type manifold due to $L$ becomes a semisymmetric  (resp. projectively, conharmonically, concircularly, conformally, and Ricci semisymmetric) manifold.
Again, the $W$-pseudosymmetric type manifold due to $L$ becomes a Deszcz pseudosymmetric (or a Ricci generalized pseudosymmetric) manifold if $L=R$, $W=R$, and $\beta=g$ (respectively, $S$).  Furthermore, if we substitute $L$ = $cir(R)$, $C$, $P$, and $har(R)$, then we get plenty of pseudosymmetric type curvature conditions. We mention that any semisymmetric manifold is also pseudosymmetric, but not vice versa.  Morris-Thorne spacetime \cite{ECS22}, Reissner-Nordstr\"{o}m  spacetime \cite{Kowa06}, Schwarzschild spacetime \cite{GP09} and Bardeen spacetime \cite{SHS_Bardeen}, for examples, are pseudosymmetric but not semisymmetric since $\mathscr{F}_L \neq 0$, but pp-wave spacetime \cite{SBK21} is semisymmetric.

\begin{defi} $($\cite{DGHZ16, DGJZ-2016, DGP-TV-2011, S09, SKH11, SK19,SYH09}$)$ 
If the rank of $(S-\varphi g)$, for a scalar $\varphi$, is $1$ (resp., $0$ and $2$), then the manifold $\mathscr{V}_n$ is said to be quasi-Einstein (resp. Einstein and $2$-quasi-Einstein). In particular, the quasi-Einstein manifold becomes Ricci simple when $\varphi =0$.
\end{defi}
	For instance,  Sultana-Dyer spacetime \cite{EDS_sultana_2022} and Som-Raychaudhuri spacetime \cite{SK16srs} are $2$-quasi-Einstein, Robertson-Walker spacetime \cite{ ARS95, ADEHM14,  DK99} and Siklos spacetime \cite{SDKC19} are quasi-Einstein; Kaigorodov spacetime \cite{SDKC19} is Einstein but Vaidya spacetime \cite{SKS19} is a Ricci simple manifold.

	


\begin{defi} $($\cite{Bess87, SK14, SK19}$)$
A semi-Riemannian manifold $\mathscr{V}_n$ is said to be Ein$(k)$ if $g$, $S$, $S^2$, $S^3$, $\cdots$, $S^k$ are linearly dependent and  $S^i$ is defined by $S^i(X,Y)=S^{i-1}(JX,Y)$ for  $i=2, 3,\cdots , k$  and $J$ be the Ricci opeartor of type $(1,1)$.
\end{defi}	
For instance, Bardeen spacetime \cite{SHS_Bardeen} and Nariai spacetime \cite{SAAC20N} are  Ein$(2)$ while Vaidya-Bonner spacetime \cite{SDC} and Lifshitz spacetime \cite{SSC19} are  Ein$(3)$ manifold.

\begin{defi}
If  a semi-Riemannian manifold $\mathscr{V}_n$ possesses the curvature relation
	\beb
	R_{efst}&=&\mathscr{R}_{11}(g\wedge g)_{efst}+\mathscr{R}_{12}(g\wedge S)_{efst}+\mathscr{R}_{22}(S \wedge S)_{efst}\\ &+& \mathscr{R}_{13}(g\wedge S^2)_{efst}+\mathscr{R}_{23}(S \wedge S^2)_{efst}+\mathscr{R}_{33}(S^2 \wedge S^2)_{efst},
	\eeb
then it is referred as a generalized Roter type manifold \cite{Desz03, DGJPZ13, DGJZ-2016, DGP-TV-2015, SK16, SK19} and it is a Roter type manifold   \cite{Desz03, DG02, DGP-TV-2011, DPSch-2013, Glog-2007} for $\mathscr{R}_{13}=\mathscr{R}_{23}=\mathscr{R}_{33}=0$.
\end{defi}
For instance, Lemaitre-Tolman-Bondi spacetime \cite{SAACD_LTB_2022} and Lifshitz spacetime \cite{SSC19} are of generalized Roter type and Melvin magnetic spacetime \cite{SAAC20} as well as Hayward spacetime \cite{SHDS_hayward}  are of Roter type.

\begin{defi} \cite{TB89, TB93}
	The manifold $(\mathscr{V}_n,g)$, is called weakly symmetric if 
	$$\nabla_d R_{efst}=\Pi_d R_{efst}+X_e R_{dfst}+X_fR_{dest}+Y_s R_{deft}+Y_t R_{defs},$$ 
where $\Pi$, $X$, $Y$ are the $1$-forms on $\mathscr{V}_n$. For $\frac{\Pi_e}{2}=X_e=Y_e$, it reduces to a Chaki pseudosymmetric manifold \cite{Chak87, Chak88}. Furthermore, if $X_e=Y_e=0$, then $\mathscr{V}_n$ turns into a recurrent manifold
 \cite{Patt52, Ruse46, Ruse49a, Ruse49b, Walk50}.	
\end{defi}

\begin{defi}  $($\cite{DD91, DGJPZ13, MM12a, MM22a, MM22b,  MM12b, MM13}$)$
	Let $\Gamma$ be a tensor of type $(0,4)$ on $\mathscr{V}_n$. If the Ricci tensor $S$ of $\mathscr{V}_n$ satisfies the relation $$\mathop{\mathcal{S}}_{e,f,s}S^d_e \otimes \Gamma_{fstd} =0,$$ then it is called a $\Gamma$-compatible and $\mathcal{S}$ is the cyclic sum over $e,f,s$.
	Similarly, if $\Gamma=C$ (resp., $R$, $har(R)$, $P$, $cir(R)$) then the Ricci tensor is conformally (resp., Riemann, conharmonic, projective, and concircular) compatible.

	Again, if the condition $$\mathop{\mathcal{S}}_{e,f,s} \nabla_e \Gamma_{fstd}=\mathop{\mathcal{S}}_{e,f,s} \Pi_e \Gamma_{fst} $$ holds, then the curvature $2$-forms for $\Gamma$ are recurrent \cite{LR89, MS12a, MS13a, MS14, Patt52}  and the $1$-forms are recurrent \cite{SKP03} if $$\nabla_eH_{fs}-\nabla_f H_{es}=\Pi_e H_{fs}-\Pi_f H_{es},$$ where $H$ is a tensor of type $(0,2)$ and $\Pi$ is a covector.
\end{defi}
For instance, the Riemann, projective, concircular, conharmonic and projective curvature tensors are compatible in Vaidya spacetime \cite{SKS19}, G\"odel spacetime \cite{DHJKS14} and Nariai spacetime \cite{SAAC20N}. Also, the conformal $2$-forms are recurrent in Hayward spacetime \cite{SHDS_hayward} and Morris-Thorne spacetime \cite{ECS22}.

\begin{defi}
	The Codazzi type \cite{F81, S81} (resp., cyclic parallel  \cite{ Gray78, SB08}) Ricci tensor of $\mathscr{V}_n$ is defined as 
	$$\nabla_e S_{fs}-\nabla_f S_{es}=0$$  $$(resp., \mathop{\mathcal{S}}_{e,f,s} \nabla_e S_{fs}=0).$$ 
\end{defi}

We note that the Ricci tensor of $(t-z)$ type plane wave spacetime \cite{EC21} is of Codazzi type and the Ricci tensor of G\"{o}del spacetime \cite{DHJKS14} is cyclic parallel.

\begin{defi} (\cite{P95, SK16, Venz85})
	Let $\Gamma$ be a tensor of type $(0,4)$ on $\mathscr{V}_n$ such that it admits the relation $$\mathop{\mathcal{S}}_{e,f,s} \Pi_e \otimes \Gamma _{fstd}=0,$$  

	where $\mathcal{L}(M)$ is the vector space of all $1$-forms with dimension greater than $1$, and $\mathcal{S}$ is the cyclic sum over $e,f,s$.  Then Venzi refers to $\mathscr{V}_n$ as a $\Gamma$-space.
	
\end{defi}
 
 For instance, generalized pp-wave spacetime \cite{SBK21} and $(t-z)$ type plane wave spacetime \cite{EC21} are Venzi space.

\begin{defi}
A semi-Riemannian manifold $\mathscr{V}_n$ admits motion about a certain vector field $\xi$ if $\pounds_{\xi} g=0$. Killing is another name for the vector field $\xi$.
\end{defi}

In $1969$, Katzin {et al.} \cite{KLD1969, KLD1970} established the idea of curvature collineation for the $(1,3)$- type Riemann curvature tensor by vanishing the Lie derivative of the tensor along a vector field. Duggal \cite{Duggal1992} again generalizes the idea of curvature collineation in $1992$ by proposing the concept of curvature inheritance for the (1,3)-type curvature tensor.

\begin{defi}\label{def_CI} (\cite{Duggal1992})
A semi-Riemannian manifold $\mathscr{V}_n$ admits curvature inheritance if it possesses the relation $$\pounds_\xi R=\zeta_1 R$$ for some vector field $\xi$, and  $\zeta_1$ is the scalar function.
Further, if $\zeta_1=0$, it becomes a curvature collineation ($\pounds_\xi R=0$).
\end{defi}
\begin{defi} \label{def_RI}(\cite{Duggal1992})
A semi-Riemannian manifold $\mathscr{V}_n$ admits Ricci inheritance if its Ricci tensor $S$ fulfills the relation $$\pounds_\xi S=\zeta_1 S$$ for some vector field $\xi$, and $\zeta_1$ is the scalar function.
Further, if $\zeta_1=0$, it becomes a Ricci collineation ($\pounds_\xi S=0$).
\end{defi}
Recently, Shaikh and Datta \cite{ShaikhDatta2022} developed the idea of curvature inheritance for $(0,4)$-type curvature tensor $R$ and presented the idea of generalized curvature inheritance, which is expressed as follows:
\begin{defi}\label{def_GCI} (\cite{ShaikhDatta2022})
If $R$ is the $(0,4)$ type curvature tensor and there is a non-Killing vector field $\xi$ that fulfills the relation $$\pounds_\xi  R=\zeta_1 R + \zeta_2 g\wedge g +\zeta_3 g\wedge S +\zeta_4 S\wedge S ,$$ then it is known as a 
 generalized curvature inheritance,  where $\zeta_1,\zeta_2,\zeta_3,\zeta_4$ are the scalar functions. In particular, $\mathscr{V}_n$ admits curvature inheritance for the $(0,4)$-type curvature tensor $R$ if $\zeta_i=0$ for $i=2,3,4$. Furthermore, it becomes curvature collineation for the $(0,4)$-type curvature tensor $R$ if $\zeta_i=0$ for $i=1,2,3,4$.\\
We note that generalized curvature inheritance holds in Robinson-Trautman spacetime \cite{ShaikhDatta2022}.
\end{defi}

\section{\bf Vaidya-Bonner-de Sitter spacetime admitting symmetry and pseudosymmetry structures} 
%
In terms of coordinates system $(t,r,\theta,\phi)$, the metric tensor  of VBdS spacetime is given as follows:
$$
g=\left(
\begin{array}{cccc}
	(1- \frac{2m(t)}{r}+\frac{q^2(t)}{r^2}-\frac{\lambda r^2}{3}) & -1 & 0 & 0 \\
	-1 & 0 & 0 & 0 \\
	0 & 0 & -r^2 & 0 \\
	0 & 0 & 0 & -r^2 \sin^2\theta
\end{array}
\right).
$$
Now, the components of the metric $g$ are  
$$\begin{array}  {c}g_{11}=\left(1- \frac{2m(t)}{r}+\frac{q^2(t)}{r^2}-\frac{\lambda r^2}{3}\right),\; g_{12}=g_{21}=-1, \\ g_{33}=-r^2, \;g_{44}=-r^2\sin^2\theta, \; g_{ij}=0, \;\mbox{otherwise}.
\end{array}$$

\indent After a rigorous computation of the components of the various curvature tensors of the VBdS spacetime, we present only the non-zero components by considering the symmetric properties of the tensors throughout the paper.\\

 The non-vanishing components of Christoffel symbols $(\Gamma^h_{ij})$ of the 2nd kind are obtained as follows:
 
\begin{eqnarray}
\begin{cases}
	\Gamma^1_{11}=-\frac{\lambda_2}{3r^3}=-\Gamma^2_{12}, \\ 
	\Gamma^2_{11}=\frac{2r^6\lambda(-3+r^2\lambda)-6(6r^2m^2+3q^2)(r^2+q^2)+m\lambda_6+3r^4m'+9r^3(q^2)'}{18r^5}, \\
	\Gamma^3_{23}=\frac{1}{r}=\Gamma^4_{24}, \ \ 
	\Gamma^1_{33}=-r, \ \
	\Gamma^2_{33}= -r+\frac{\lambda_1}{3r}=\frac{1}{\sin^2\theta}\Gamma^2_{44}, \ \ 
	\Gamma^4_{34}=\cot\theta, \\
	\Gamma^1_{44}=-r \sin^2 \theta, \ \
	
	 	\Gamma^3_{44}=-\cos\theta \sin\theta. \\
	 	\end{cases}
\end{eqnarray}
The non-vanishing components of the Riemann and Ricci curvature tensors are obtained as follows:
\begin{eqnarray}
\begin{cases}
	R_{1212}=\frac{\lambda_3}{3r^4},
	 \ \ 
	
	R_{1313}=-\frac{(-3r^2
		+\lambda_1)\lambda_2}{9r^4}-m'+\frac{(q^2)'}{2r}=\frac{1}{\sin^2\theta}R_{1414}, \\
	
	R_{1323}=-\frac{\lambda_2}{3r^2}=\frac{1}{\sin^2\theta}R_{1424}, \ \

	R_{3434}=-\frac{\lambda_1}{3} \sin^2 \theta;\\
\end{cases}
\end{eqnarray}
\begin{eqnarray}
\begin{cases}
	S_{11}=\frac{r^6 \lambda(3-r^2\lambda)-\lambda_4-6r^4m'+3r^3(q^2)'}{3r^6}, \\ 
	S_{22}= -\lambda+\frac{q^2}{r^4}, \ \
	S_{33}=-\frac{r^4\lambda+q^2}{r^2}, \ \ 
	S_{44}=\sin^2\theta \ S_{33};
\end{cases}
\end{eqnarray}
where $\lambda_1= r^4\lambda+6rm-3q^2$, $\lambda_2=r^4\lambda-3rm+3q^2$, $\lambda_3=r^4 \lambda+6rm-9q^2$, $\lambda_4=3r^2q^2-4r^4\lambda q^2+3q^4+6r^5m\lambda-6mrq^2$, $\lambda_5=r^4\lambda-12rm+9q^2$,  $\lambda_6=3r^3+r^5\lambda+9rq^2$, $\lambda_7=3r^2-26r^4 \lambda+3q^2$ and $\lambda_8=3r^2-r^4\lambda+3q^2$ with $m'=\frac{d}{dt}(m(t))$ and $q'=\frac{d}{dt}(q(t))$.\\
Also,  the scalar curvature $\bar{\kappa}$ is given by, $\bar{\kappa}= 4\lambda.$\\ 
From the above calculation, we can state the following:

\begin{pr}\label{pr1}
	The VBdS spacetime is neither an Einstein manifold nor a quasi-Einstein manifold, but rather 
	$(i)$ it is $2$-quasi-Einstein as rank $(S-\varphi g)=2$ for $\varphi=\frac{r^4 \lambda+q^2}{r^4}$. \\
	$(ii)$ In the VBdS spacetime the common form of a tensor compatible with $R$ is  \\

$$	\left(
	\begin{array}{cccc}
	\mathcal{H}_{11} & \mathcal{H}_{12} & 0 & 0 \\
	\mathcal{H}_{12}-\frac{3r\mathcal{H}_{22}(2rm'-(q^2)')}{2\lambda_2} & \mathcal{H}_{22} & 0 & 0 \\
	0 & 0 & \mathcal{H}_{33} & \mathcal{H}_{34} \\
	0 & 0 & \mathcal{H}_{34} & \mathcal{H}_{44}
	\end{array}
	\right)
	$$

where $\mathcal{H}_{ij}$ are arbitrary  scalars. 	
\end{pr}

The non-vanishing components of the second rank Ricci tensor $S^2_{ij}$  are determined as follows:
\begin{eqnarray} \label{S}
\begin{cases}
 S^2_{11}=-\frac{(r^4\lambda-q^2)\lambda_4+(r^8\lambda^2+12r^4m'-3r^6\lambda-6r^3(q^2)')}{3r^{10}}, \\
S^2_{12}=-(\lambda-\frac{q^2}{r^4})^2, \ \
S^2_{33}=-\frac{(r^4\lambda+q^2)^2}{r^6}, \\
S^2_{44}=\sin^2\theta \ S^2_{33}.
\end{cases}
\end{eqnarray}
\indent Let $\mathscr{W}^1=g\wedge g$, $\mathscr{W}^2=g\wedge S$, $\mathscr{W}^3=S \wedge S $, $\mathscr{W}^4=g\wedge S^2$, $\mathscr{W}^5=S \wedge S^2$  and $\mathscr{W}^6=S^2\wedge S^2$. Then the Kulkarni-Nomizu product tensors $\mathscr{W}^1$, $\mathscr{W}^2$, $\mathscr{W}^3$, $\mathscr{W}^4$, $\mathscr{W}^5$ and $\mathscr{W}^6$ are computed as below:

\begin{eqnarray}
\begin{cases} \label{1}
\mathscr{W}^1_{1212}=2, \ \
\mathscr{W}^1_{1313}=2(r^2-\frac{\lambda_1}{3})=\frac{1}{\sin^2\theta}\mathscr{W}^1_{1414}, \\
\mathscr{W}^1_{1323}=-2r^2, \ \

\mathscr{W}^1_{1424}= \sin^2\theta \  \mathscr{W}^1_{1323}, \ \
\mathscr{W}^1_{3434}=r^2 \  \mathscr{W}^1_{1424}; \\
\end{cases}
\end{eqnarray}

\begin{eqnarray}
\begin{cases}   \label{3}
\mathscr{W}^2_{1212}=2\lambda-\frac{2q^2}{r^4}, \\
\mathscr{W}^2_{1313}=2r^2\lambda-\frac{2r^4\lambda^2}{3}-4r\lambda m+2\lambda q^2-2m'+\frac{(q^2)'}{r}=\frac{1}{\sin^2\theta}\mathscr{W}^2_{1414}, \\
\mathscr{W}^2_{1323}=-2r^2 \lambda =\frac{1}{\sin^2\theta}\mathscr{W}^2_{1424}, \ \
\mathscr{W}^2_{3434}=-2(r^4 \lambda+q^2) \sin^2 \theta; \\
\end{cases}
\end{eqnarray}

\begin{eqnarray}
\begin{cases}  \label{2}
\mathscr{W}^3_{1212}=2(\lambda - \frac{q^2}{r^4})^2, \\
\mathscr{W}^3_{1313}=\frac{2(r^4\lambda+q^2)(r^6\lambda(3-r^2\lambda)-\lambda_4-6r^4m'+3(q^2)')}{3r^8}=\frac{1}{\sin^2\theta}\mathscr{W}^3_{1414},  \\
\mathscr{W}^3_{1323}=\frac{2(-r^8\lambda^2+q^4)}{r^6}=\frac{1}{\sin^2\theta}\mathscr{W}^3_{1424}, \\ \ \ 
 
\mathscr{W}^3_{3434}=-\frac{2(r^4\lambda+q^2)^2 \sin^2 \theta}{r^4}; \\
\end{cases}
\end{eqnarray}

\begin{eqnarray}
\begin{cases}   \label{4}
\mathscr{W}^4_{1212}=2(\lambda-\frac{q^2}{r^4})^2, \\
\mathscr{W}^4_{1313}=\frac{-2(-3r^2+\lambda_1)(r^8\lambda^2+q^4)+12r^4(q^2-r^4\lambda)m'+6r^3(r^4\lambda-q^2)(q^2)'}{3r^8}=\frac{1}{\sin^2\theta}\mathscr{W}^4_{1414}, \\
\mathscr{W}^4_{1323}=-\frac{2(r^8\lambda^2+q^4)}{r^6}=\frac{1}{\sin^2\theta}\mathscr{W}^4_{1424}, \\

\mathscr{W}^4_{3434}=-\frac{2(r^4\lambda+q^2)^2 \sin^2 \theta}{r^4}; \\
\end{cases}
\end{eqnarray}

\begin{eqnarray}
\begin{cases}   \label{5}
\mathscr{W}^5_{1212}=\frac{2(r^4\lambda-q^2)^3}{r^{12}}, \\
\mathscr{W}^5_{1313}=-\frac{(r^4\lambda+q^2)(2r\lambda)(-3r^2+\lambda_1)(r^4\lambda-q^2)+6r(3r^4\lambda-q^2)m'+3(-3r^4\lambda+q^2)(q^2)'}{3r^9}=\frac{1}{\sin^2\theta}\mathscr{W}^5_{1414}, \\
\mathscr{W}^5_{1323}=-2r^2 \lambda^3+\frac{2\lambda q^4}{r^6}=\frac{1}{\sin^2\theta}\mathscr{W}^5_{1424}, \\

\mathscr{W}^5_{3434}=-\frac{2(r^4 \lambda +q^2)^3 \sin^2 \theta}{r^8}; \\

\end{cases}
\end{eqnarray}

\begin{eqnarray}
\begin{cases}   \label{6}
\mathscr{W}^6_{1212}=2(\lambda- \frac{q^2}{r^4})^4, \\
\mathscr{W}^6_{1313}=-\frac{2(r^4\lambda -q^2)(r^4\lambda+q^2)^2 \lambda_4+r^8\lambda^2-3r^6\lambda+12r^4m'-6r^3(q^2)'}{3r^{16}}=\frac{1}{\sin^2\theta}\mathscr{W}^6_{1414},  \\
\mathscr{W}^6_{1323}=-\frac{2(-r^8\lambda^2+q^4)^2}{r^{14}}=\frac{1}{\sin^2\theta}\mathscr{W}^6_{1424}, \\

\mathscr{W}^6_{3434}=-\frac{2(r^4 \lambda+q^2)^4 \sin^2 \theta}{r^{12}}. \\

\end{cases}
\end{eqnarray}
\indent In view of  $\eqref{1} - \eqref{6}$, we have deduced the following relation of linear dependency stated as follows:

\begin{eqnarray*}
R_{efst}&=&\mathscr{R}_{11}(g\wedge g)_{efst}+\mathscr{R}_{12}(g\wedge S)_{efst}+\mathscr{R}_{22}(S \wedge S)_{efst}\nonumber\\ &+&
\mathscr{R}_{13}(g\wedge S^2)_{efst}+\mathscr{R}_{23}(S \wedge S^2)_{efst}+\mathscr{R}_{33}(S^2\wedge S^2)_{efst},
\end{eqnarray*}

	where $\mathscr{R}_{11}$, $\mathscr{R}_{12}$, $\mathscr{R}_{22}$, $\mathscr{R}_{13}$, $\mathscr{R}_{23}$ and $\mathscr{R}_{33}$ are computed as given below:
$$\begin{array}{c}	
\mathscr{R}_{11}= 1,
\mathscr{R}_{12}=1, \\
\mathscr{R}_{22}=\frac{r^4(9r^{17}\lambda^4m-9r^{13}\lambda^4q^2-2r^9\lambda^2(2r^3(12+7\lambda)-15m)q^4-30r^8\lambda^2q^6+3r(4r^3\lambda+3m)q^8-9q^{10})}{12(q^6-r^8q^2\lambda^2)^2},\\

\mathscr{R}_{13}=\frac{r^4(9r^9\lambda^2m-9r^5\lambda(r^3\lambda+2m)q^2-r(4r^3(6-\lambda)+3m)q^4-3q^6)}{12q^4(r^4\lambda+q^2)^2}, \\

\mathscr{R}_{23}=\frac{r^8(-18r^{17}\lambda^4m+18r^{13}\lambda^3(r^3\lambda+m)q^2+r^9\lambda^2(r^3(48+7\lambda)-30m)q^4+3r^5\lambda(11r^3\lambda-6m)q^6+9r^4\lambda q^8-3q^{10})}{6q^4(-r^4\lambda+q^2)^2(r^4\lambda+q^2)^3}, \\

\mathscr{R}_{33}=\frac{r^{13}(9r^{12}\lambda^3m-9r^8\lambda^2(r^3\lambda+m)q^2+9r^4\lambda(-2r^3+m)q^4+r^3(6-7\lambda+3m)q^6)}{6q^4(-r^4\lambda+q^2)^2(r^4\lambda+q^2)^3}.
\end{array}$$
Contracting the last relation, we get the following:

\begin{eqnarray*}
S^3+\frac{q^2-3r^4 \lambda}{r^4} S^2+\frac{(3r^4 \lambda+q^2)(r^4\lambda-q^2)}{r^8} S+ \frac{-(r^4 \lambda-q^2)^2 (r^4\lambda+q^2)}{r^{12}} g=0.
\end{eqnarray*}
This leads to the following:
\begin{pr}  \label{pr2}
The nature of the VBdS spacetime is neither an Einstein manifold of level $2$ nor Roter type, but rather  $(i)$ an Einstein manifold of level $3$ and $(ii)$ generalized Roter type.   
\end{pr}
Let $\mathscr{P}^1=C_{efst}$ and $\mathscr{P}^2=\nabla C_{efst,d}$.
Then the tensors $\mathscr{P}^1$ and $\mathscr{P}^2$ are obtained as follows:
\begin{eqnarray}
\begin{cases}
	\mathscr{P}^1_{1212}=\frac{2rm-2q^2}{r^4}, \ \ 
	\mathscr{P}^1_{1313}=\frac{(-3r^2+\lambda_1)(rm-q^2)}{3r^4}=\frac{1}{\sin^2 \theta}\mathscr{P}^1_{1414}, \\
	\mathscr{P}^1_{1323}=\frac{rm-q^2}{r^2}=\frac{1}{\sin^2 \theta}\mathscr{P}^1_{1424}, \ \
	
	\mathscr{P}^1_{3434}=2(-rm+q^2) \ \sin^2\theta;\\
\end{cases}
\end{eqnarray}

\begin{eqnarray}
\begin{cases}
	 \mathscr{P}^2_{1212,1}= \frac{2rm'-2(q^2)'}{r^4}, \ \ 
	 \mathscr{P}^2_{1212,2}=\frac{-6rm+8q^2}{r^5}, \ \
 \mathscr{P}^2_{1213,3}=-\frac{(-3r^2+\lambda_1)(rm-q^2)}{r^5}=\frac{1}{\sin^2 \theta}\mathscr{P}^2_{1214,4}, \\   
 \mathscr{P}^2_{1223,3}=\frac{-3rm+3q^2}{r^3}=\frac{1}{\sin^2 \theta}\mathscr{P}^2_{1224,4},\ \   
 \mathscr{P}^2_{1313,1}=\frac{(-3r^2+\lambda_1)(rm'-(q^2)')}{3r^4}=\frac{1}{\sin^2 \theta}\mathscr{P}^2_{1414,1}, \\
\mathscr{P}^2_{1313,2}=-\frac{(3rm-4q^2)(-3r^2+\lambda_1)}{3r^5}=\frac{1}{\sin^2 \theta}\mathscr{P}^2_{1414,2}, \ \	
	\mathscr{P}^2_{1323,1}=\frac{rm'-(q^2)'}{r^2}=\frac{1}{\sin^2 \theta}\mathscr{P}^2_{1424,1}, \\  \mathscr{P}^2_{1323,2}=\frac{-3rm+4q^2}{r^3}=\frac{1}{\sin^2 \theta}\mathscr{P}^2_{1424,2}, \ \
	\mathscr{P}^2_{2334,4}=\frac{3(-rm+q^2)\sin^2\theta}{r}=\mathscr{P}^2_{2434,3}, \\
		
		\mathscr{P}^2_{3434,1}=-2\sin^2\theta (rm'-(q^2)'), \ \
		\mathscr{P}^2_{3434,2}=\frac{2(3rm-4q^2)\sin^2\theta}{r}.    
\end{cases}
\end{eqnarray}

The above calculation yields the following: 

\begin{pr} \label{pr3}
The VBdS spacetime with $q^2-rm\neq 0$ is not  conformally recurrent but $(i)$ its curvature $2$-forms are recurrent for  $\Gamma=\left\lbrace \frac{rm'-(q^2)'}{rm-q^2}, \frac{q^2}{r^2m-rq^2},0,0\right\rbrace $, and $(ii)$  the common form of a tensor compatible with $C$ is  \\

$$	\left(
\begin{array}{cccc}
\mathcal{H}_{11} & \mathcal{H}_{12} & 0 & 0 \\
\mathcal{H}_{12} & \mathcal{H}_{22} & 0 & 0 \\
0 & 0 & \mathcal{H}_{33} & \mathcal{H}_{34} \\
0 & 0 & \mathcal{H}_{34} & \mathcal{H}_{44}
\end{array}
\right)
$$
where $\mathcal{H}_{ij}$ are arbitrary  scalars.
\end{pr}

Let $\mathscr{P}^3=cir(R)_{efst}$, $\mathscr{P}^4= har(R)_{efst}$ and $\mathscr{P}^5=P_{efst}$.
Then the tensors $\mathscr{P}^3$, $\mathscr{P}^4$ and $\mathscr{P}^5$ are calculated as follows:

\begin{eqnarray}
\begin{cases} \label {W}
\mathscr{P}^3_{1212}=\frac{2rm-3q^2}{r^4}, \\

\mathscr{P}^3_{1313}=\frac{12r^2m^2+(6r^2-2r^4\lambda)q^2+6q^4-2m\lambda_6+3r^3(-2rm'+(q^2)')}{6r^4}=\frac{1}{\sin^2\theta}\mathscr{P}^3_{1414}, \\
\mathscr{P}^3_{1323}=\frac{rm-q^2}{r^2}=\frac{1}{\sin^2\theta}\mathscr{P}^3_{1424}, \ \
\mathscr{P}^3_{3434}=(-2rm+q^2) \sin^2 \theta; \\
\end{cases}
\end{eqnarray}

\begin{eqnarray}
\begin{cases}  \label {K}
\mathscr{P}^4_{1212}=-\frac{2 \lambda}{3}+\frac{2rm-2q^2}{r^4}, \\
\mathscr{P}^4_{1313}=\frac{(2r^4\lambda+3rm-3q^2)(-3r^2+\lambda_1)}{9r^4}=\frac{1}{\sin^2\theta}\mathscr{P}^4_{1414}, \\
\mathscr{P}^4_{1323}=\frac{2r^4 \lambda+3rm-3q^2}{3r^2}=\frac{1}{\sin^2\theta}\mathscr{P}^4_{1424}, \ \
\mathscr{P}^4_{3434}=\frac{2}{3} \lambda_2 \sin^2 \theta; 
\end{cases}
\end{eqnarray}

\begin{eqnarray}
\begin{cases}  \label {P}
\mathscr{P}^5_{1211}=\frac{2rm'-(q^2)')}{3r^3}, \ \
\mathscr{P}^5_{1212}=\frac{6rm-8q^2}{3r^4}=\mathscr{P}^5_{1221}, \\

\mathscr{P}^5_{1313}=\frac{36r^2m^2-8r^2(-3+r^2\lambda)q^2+24q^4+6m(-3r^3+r^5\lambda-11rq^2)+3r^3(-2rm'+(q^2)')}{18r^4}=\frac{1}{\sin^2\theta}\mathscr{P}^5_{1414},\\

\mathscr{P}^5_{1323}=\frac{3rm-4q^2}{3r^2}=\frac{1}{\sin^2\theta}\mathscr{P}^5_{1424}=\mathscr{P}^5_{2313}=\frac{1}{\sin^2\theta}\mathscr{P}^5_{2414}, \\ \ 
\mathscr{P}^5_{1331}=\frac{-36r^2m^2+4r^2(-3+r^2\lambda)q^2-12q^4+6rm(3r^2-r^4\lambda+7q^2)+9r^3(2rm'-(q^2)')}{18r^4}=\frac{1}{\sin^2\theta}\mathscr{P}^5_{1441},\\

\mathscr{P}^5_{1332}=\frac{-3rm+2q^2}{3r^2}=\frac{1}{\sin^2\theta}\mathscr{P}^5_{1442}=\mathscr{P}^5_{2331}=\frac{1}{\sin^2\theta}\mathscr{P}^5_{2441},\\ \ \
\mathscr{P}^5_{3434}=\frac{2}{3} (-3rm+2q^2) \sin^2 \theta=-\mathscr{P}^5_{3443}. \\

\end{cases}
\end{eqnarray}

From the equation $\eqref{W}$ - $\eqref{P}$, we have the following: 
\begin{pr} \label{pr4}
	The VBdS spacetime fulfills the following curvature conditions:
	 \begin{enumerate}[label=(\roman*)]
	 	\item  The common form of a tensor compatible with $cir(R)$ is 	
	 	$$	\left(
	 	\begin{array}{cccc}
	 		\mathcal{H}_{11} & \mathcal{H}_{12} & 0 & 0 \\
	 		\mathcal{H}_{12}+\frac{r\mathcal{H}_{22}(2rm'-(q^2)')}{2(rm-q^2)} & \mathcal{H}_{22} & 0 & 0 \\
	 		0 & 0 & \mathcal{H}_{33} & \mathcal{H}_{34} \\
	 		0 & 0 & \mathcal{H}_{34} & \mathcal{H}_{44}
	 	\end{array}
	 	\right)
	 	$$
        
         \item The common form of a tensor compatible with $har(R)$ is 	$$	\left(
        \begin{array}{cccc}
        \mathcal{H}_{11} & \mathcal{H}_{12} & 0 & 0 \\
        \mathcal{H}_{12} & \mathcal{H}_{22} & 0 & 0 \\
        0 & 0 & \mathcal{H}_{33} & \mathcal{H}_{34} \\
        0 & 0 & \mathcal{H}_{34} & \mathcal{H}_{44}
        \end{array}
        \right)
        $$

        \item  The common form of a tensor compatible with $P$ is 	 	 
        $$	\left(
        \begin{array}{cccc}
        \mathcal{H}_{11} & \mathcal{H}_{12} & 0 & 0 \\
        \mathcal{H}_{12}+\frac{3r\mathcal{H}_{22}(2rm'-(q^2)')}{2(3rm-2q^2)} & \mathcal{H}_{22} & 0 & 0 \\
        0 & 0 & \mathcal{H}_{33} & \mathcal{H}_{34} \\
        0 & 0 & \mathcal{H}_{34} & \mathcal{H}_{44}
        \end{array}
        \right)
        $$
    \end{enumerate}
    where $\mathcal{H}_{ij}$ are arbitrary  scalars.
\end{pr}

 Again, let $\mathscr{W}^7=R\cdot R$, $\mathscr{W}^8=C\cdot C$, $\mathcal{G}^1=Q(g,R)$, $\mathcal{G}^2=Q(S,R)$ and $\mathcal{G}^3=Q(g,C)$, then the  components of $\mathscr{W}^7$, $\mathscr{W}^8$, $\mathcal{G}^1$, $\mathcal{G}^2$ and $\mathcal{G}^3$ are obtained as given below: 

\begin{eqnarray}
\begin{cases} \label {R R}
\mathscr{W}^7_{1313,12}=-\frac{\lambda_3 (2rm'-(q^2)')}{3r^5}=\frac{1}{\sin^2\theta}\mathscr{W}^7_{1414,12}, \ \

\mathscr{W}^7_{1213,13}=\frac{\lambda_3 (2rm'-(q^2)')}{6r^5}=\frac{1}{\sin^2\theta}\mathscr{W}^7_{1214,14},\\
\mathscr{W}^7_{1223,13}=-\frac{\lambda_2 (3rm-4q^2)}{3r^6}=-\mathscr{W}^7_{1213,23}=\frac{1}{\sin^2\theta}\mathscr{W}^7_{1224,14}=\frac{1}{\sin^2\theta}\mathscr{W}^7_{1214,24}, \\
\mathscr{W}^7_{1434,13}=\frac{\sin^2 \theta(2(-3r^2+\lambda_1)(3rm-2q^2)\lambda_2-6r^4 \lambda_4+m'+3r^3\lambda_4(q^2)')}{18r^6}=-\mathscr{W}^7_{1334,14} \\
\mathscr{W}^7_{2434,13}=\frac{(3rm-2q^2)\lambda_2 \sin^2\theta}{3r^4} =-\mathscr{W}^7_{2334,14}=\mathscr{W}^7_{1434,23}=-\mathscr{W}^7_{1334,24}; \\

\end{cases}
\end{eqnarray}

\begin{eqnarray}
\begin{cases}  \label {C C}
\mathscr{W}^8_{1223,13}=\frac{3(-rm+q^2)^2}{r^6}=-\mathscr{W}^8_{1213,23}, \\ 
\mathscr{W}^8_{1434,13}=-\frac{(-3r^2+\lambda_1)(-rm+q^2)^2 \sin^2\theta}{r^6}=-\mathscr{W}^8_{1334,14}, \\   
\mathscr{W}^8_{2434,13}=-\frac{3(-rm+q^2)^2 \sin^2 \theta}{r^4}=-\mathscr{W}^8_{2334,14}=\mathscr{W}^8_{1434,23}=-\mathscr{W}^8_{1334,24},    \\

\mathscr{W}^8_{1224,14}=\frac{3(-rm+q^2)^2 \sin^2 \theta}{r^6}=-\mathscr{W}^8_{1214,24}; \\

\end{cases}
\end{eqnarray}

\begin{eqnarray}
\begin{cases}   \label {Q(g,R)}	
\mathcal{G}^1_{1313,12}=-2m'+\frac{(q^2)'}{r}=\frac{1}{\sin^2\theta}\mathcal{G}^1_{1414,12}, \ \

\mathcal{G}^1_{1213,13}=m'-\frac{(q^2)'}{2r}=\frac{1}{\sin^2\theta}\mathcal{G}^1_{1214,14}, \\
\mathcal{G}^1_{1223,13}=-\frac{-3rm+4q^2}{r^2}=-\mathcal{G}^1_{1213,23}, \ \
\mathcal{G}^1_{1214,24}=\frac{(3rm-4q^2)\sin^2\theta}{r^2}=- \mathcal{G}^1_{1224,14}, \\
\mathcal{G}^1_{1434,13}=\frac{\sin^2 \theta (36r^2m^2-4r^2(-3+r^2\lambda)q^2+12q^4+6m(\lambda_6-16rq^2)+3r^3(-2rm'+(q^2)')}{6r^2}-\mathcal{G}^1_{1334,14}, \\
\mathcal{G}^1_{2434,13}=(3rm-2q^2) \sin^2 \theta=-\mathcal{G}^1_{2334,14}=\mathcal{G}^1_{1434,23}=-\mathcal{G}^1_{1334,24}; \\

\end{cases}
\end{eqnarray}

\begin{eqnarray}
\begin{cases}   \label {Q(S,R)}
\mathcal{G}^2_{1313,12}=-\frac{ \lambda_3 (2rm'-(q^2)')}{3r^5}=\frac{1}{\sin^2\theta}\mathcal{G}^2_{1414,12}, \ \

\mathcal{G}^2_{1213,13}=\frac{ \lambda_3 (2rm'-(q^2)')}{6r^5}=\frac{1}{\sin^2\theta}\mathcal{G}^2_{1214,14}, \\
\mathcal{G}^2_{1223,13}=\frac{-3rm(3r^4\lambda+q^2)+2q^2(5r^4\lambda+3q^2)}{3r^6}=\frac{1}{\sin^2\theta}\mathcal{G}^2_{1224,14}=\mathcal{G}^2_{1213,23}=\frac{1}{\sin^2\theta}\mathcal{G}^2_{1214,24}, \\
\mathcal{G}^2_{1434,13}=\frac{\sin^2 \theta(2(18m^2(3r^5\lambda-rq^2)+3m(3r^6\lambda(r^2\lambda-3)+q^2\lambda_7+r^3(8\lambda q^2\lambda_8-3(r^4\lambda+9q^2)m'))+3r^2 \lambda_5 (q^2)'))}{18r^5}=-\mathcal{G}^2_{1334,14},  \\
\mathcal{G}^2_{2434,13}=\frac{(9r^4\lambda m-(8r^3\lambda +3m)q^2)\sin^2\theta}{3r^3}=\mathcal{G}^2_{1434,23}, \\

\mathcal{G}^2_{2334,14}=\frac{1}{3}(8\lambda q^2+m(-9r\lambda+\frac{3q^2}{r^3}))\sin^2 \theta=\mathcal{G}^2_{1334,24}; \\
\end{cases}
\end{eqnarray}

\begin{eqnarray}
\begin{cases}   \label {Q(g,C)}
\mathcal{G}^3_{1223,13}=\frac{-3rm+3q^2}{r^2}=-\mathcal{G}^3_{1213,23}, \\
\mathcal{G}^3_{1434,13}=\frac{(-3r^2+\lambda_1)(rm-q^2)\sin^2\theta}{r^2}=-\mathcal{G}^3_{1334,14}, \\
\mathcal{G}^3_{2434,13}=3(rm-q^2)\sin^2\theta=-\mathcal{G}^3_{2334,14}=\mathcal{G}^3_{1434,23}=-\mathcal{G}^3_{1334,24}, \\
\mathcal{G}^3_{1224,14}=\frac{3(-rm+q^2)\sin^2\theta}{r^2}=-\mathcal{G}^3_{1214,24}. \\

\end{cases}
\end{eqnarray}

By virtue of $\eqref{R R} - \eqref{Q(g,C)}$, we can assert the following:
\begin{pr} \label{pr5}
	The VBdS spacetime with $q^2-rm\neq 0$ is neither a Ricci generalized pseudosymmetric manifold nor a pseudosymmetric manifold, but 
	\begin{enumerate}[label=(\roman*)]
		\item it admits $C\cdot C=-\frac{rm-q^2}{r^4} Q(g,C)$,
		
		\item it satisfies the curvature condition
		 $$R\cdot R+\beta Q(g,C)=Q(S,R),$$ where $\beta=\frac{2r^5\lambda m+3r^2m^2-2r^4\lambda q^2-6rmq^2+2q^4}{3r^4(rm-q^2)} $.
	\end{enumerate}
		
\end{pr}
\begin{note}
If $q^2-rm=0$, then the VBdS spacetime satisfies $C\cdot C=0$.

\end{note}

 Also, if $\mathscr{W}^9$=$R\cdot C$, $\mathcal{G}^4$=$Q(S,C)$ and $\mathscr{W}^{10}$=$C\cdot R$, then the components of $\mathscr{W}^9$, $\mathcal{G}^4$ and $\mathscr{W}^{10}$ are computed as given below:

\begin{eqnarray}
\begin{cases} \label{R C}
\mathscr{W}^9_{1213,13}=\frac{3(rm-q^2)(2rm'-(q^2)')}{2r^5}=\frac{1}{\sin^2\theta}\mathscr{W}^9_{1214,14}, \\
\mathscr{W}^9_{1223,13}=-\frac{(rm-q^2)\lambda_2}{r^6}=\frac{1}{\sin^2\theta}\mathscr{W}^9_{1224,14}=\frac{1}{\sin^2\theta}\mathscr{W}^9_{1214,24}=-\mathscr{W}^9_{1213,23}, \\
\mathscr{W}^9_{1434,13}=\frac{(rm-q^2)\sin^2 \theta (r^6 \lambda (-3+r^2 \lambda)-18r^2m^2-9r^2q^2-9q^4+3m\lambda_6+9r^4m'-9r^3qq')}{3r^6}=-\mathscr{W}^9_{1334,14}, \\

\mathscr{W}^9_{2434,13}=\frac{(rm-q^2)\lambda_2 \sin^2 \theta}{r^4}=-\mathscr{W}^9_{2334,14}=\mathscr{W}^9_{1434,23}=-\mathscr{W}^9_{1334,24}; \\

\end{cases}
\end{eqnarray}

\begin{eqnarray}
\begin{cases} \label{Q(S,C)}

\mathcal{G}^4_{1313,12}=-\frac{2(rm-q^2)(2rm'-(q^2)')}{r^5}=\frac{1}{\sin^2\theta}\mathcal{G}^4_{1414,12} , \\

\mathcal{G}^4_{1213,13}=\frac{(rm-q^2)(2rm'-(q^2)')}{r^5}=\frac{1}{\sin^2\theta}\mathcal{G}^4_{1214,14} , \\ 
\mathcal{G}^4_{1223,13}=-\frac{(rm-q^2)(3r^4\lambda+q^2)}{r^6}=\frac{1}{\sin^2\theta}\mathcal{G}^4_{1224,14}=-\frac{1}{\sin^2\theta}\mathcal{G}^4_{1214,22}=-\mathcal{G}^4_{1213,23},\\
\mathcal{G}^4_{1434,13}=\frac{(rm-q^2)\sin^2 \theta (3r^6 \lambda(-3+r^2 \lambda)+r^2(3-10r^2 \lambda)q^2+3q^4+6m(3r^5 \lambda -rq^2)+12r^4m'-12r^3qq')}{3r^6}=-\mathcal{G}^4_{1334,14}, \\
\mathcal{G}^4_{2434,13}=\frac{(3r^4\lambda-q^2)(rm-q^2)\sin^2 \theta}{r^4}=-\mathcal{G}^4_{2334,14}=\mathcal{G}^4_{1434,23}=-\mathcal{G}^4_{1334,24}; \\

\end{cases}
\end{eqnarray}

\begin{eqnarray}
\begin{cases} \label{C R}
\mathscr{W}^{10}_{1313,12}=-\frac{2(rm-q^2)(2rm'-(q^2)')}{r^5}=\frac{1}{\sin^2\theta}\mathscr{W}^{10}_{1414,12}, \\

\mathscr{W}^{10}_{1213,13}=-\frac{(rm-q^2)(2rm'-(q^2)')}{2r^5}=\frac{1}{\sin^2\theta}\mathscr{W}^{10}_{1214,14}, \\
\mathscr{W}^{10}_{1223,13}=\frac{(3rm-4q^2)(rm-q^2)}{r^6}=\frac{1}{\sin^2\theta}\mathscr{W}^{10}_{1224,14}=-\frac{1}{\sin^2\theta}\mathscr{W}^{10}_{1214,24}=-\mathscr{W}^{10}_{1213,23}, \\
\mathscr{W}^{10}_{1434,13}=-\frac{(rm-q^2)\sin^2 \theta (18r^2m^2+(6r^2-2r^4 \lambda)q^2+6q^4+3m( \lambda_6 -16rq^2)-3r^4m'+3r^3qq')}{3r^6}=-\mathscr{W}^{10}_{1334,14}, \\
\mathscr{W}^{10}_{2434,13}=-\frac{(3rm-2q^2)(rm-q^2)\sin^2 \theta}{r^4}=-\mathscr{W}^{10}_{2334,14}=\mathscr{W}^{10}_{1434,23}=-\mathscr{W}^{10}_{1334,24}. \\ 

\end{cases}
\end{eqnarray}

From the equation $\eqref{R C}$ - $\eqref{C R}$, it represents that the generalized Einstein metric condition is not satisfied by the VBdS spacetime, but rather it fulfills the following:  \\
$$R\cdot C+C\cdot R=Q(S,C)+\frac{-2 (r^4 \lambda +3rm-3q^2)}{3r^4} Q(g,C).$$

From the Proposition \ref{pr1}-- Proposition \ref{pr5}, we can derive the following results about VBdS spacetime:

\begin{thm}
 The VBdS spacetime with $q^2-rm\neq0$ admits the following geometric structures:
\begin{enumerate}[label=(\roman*)]
\item it realizes the relation  $C\cdot C=-\frac{rm-q^2}{r^4} Q(g,C)$, and  hence  \\ $C\cdot har(R)=-\frac{rm-q^2}{r^4} Q(g,har(R))$,
 
\item it realizes the relation  $har(R)\cdot C=-\frac{2r^4 \lambda+3rm-3q^2}{3r^4}Q(g,C)$,  and hence \\ $har(R)\cdot har(R)=-\frac{2r^4 \lambda+3rm-3q^2}{3r^4}Q(g,har(R))$, 
\item it indulges the pseudosymmetric type relation $$R\cdot R+\beta Q(g,C)=Q(S,R),$$ where $\beta=\frac{2r^5\lambda m+3r^2m^2-2r^4\lambda q^2-6rmq^2+2q^4}{3r^4(rm-q^2)} $,	
\item its conformal $2$-forms are recurrent for the $1$-form $\Gamma=\left\lbrace \frac{rm'-(q^2)'}{rm-q^2}, \frac{q^2}{r^2m-rq^2},0,0\right\rbrace $, 
\item it is a spacetime of generalized Roter type and Einstein of level $3$, 
\item its Ricci tensor is compatible for $R$, $P$, $cir(R)$, $har(R)$ and $C$.

\end{enumerate}

\end{thm}

\begin{rem}
The VBdS spacetime does not admit the following geometric structures:
	\begin{enumerate}[label=(\roman*)]
		\item any semisymmetric type conditions,
		
		\item Deszcz's pseudosymmetric, Ricci generalized pseudosymmetric type conditions,
		
		\item  $C\cdot R= \mathscr{F}_R Q(g,R)$, for any smooth function $\mathscr{F}_R$ and consequently projectively pseudosymmetric, concircularly pseudosymmetric,
		
		\item Venzi space for $R$, $C$, $har(R)$, $cir(R)$ and $P$,
		
		\item Einstein or quasi-Einstein  condition,
		
		\item Chaki pseudosymmetry.
	\end{enumerate}	
\end{rem}

Since, VBdS spacetime is the generalization of Vaidya-Bonner spacetime, Vaidya spacetime, and Schwarzschild spacetime, as a particular case of our result, we can obtain the curvature properties of such spacetimes. For $\lambda=0$ in (\ref{VBdS}), BVdS spacetime reduces to Vaidya-Bonner spacetime and hence we can propound the following:

\begin{cor} \cite{SDC}
The Vaidya-Bonner spacetime with $q^2-rm\neq 0$ has the following curvature restricted geometric properties:
    \begin{enumerate}[label=(\roman*)]
        \item its scalar curvature vanishes, i.e.,  $\bar{\kappa}=0$, and consequently $R=cir(R)$ as well as $C=har(R)$,
        
        \item it admits   pseudosymmetric Weyl conformal tensor,

        \item it satisfies the curvature conditions $R\cdot R-\frac{3r^2m^2-6rmq^2+2q^4}{3r^4(rm-q^2)} Q(g,C)=Q(S,R)$ and \\ $R\cdot C+C\cdot R=Q(S,C)+\frac{2(rm-q^2)}{r^4}Q(g,C)$,

        \item it is $2$-quasi-Einstein  for $\varphi=-\frac{q^2}{r^4}$,  pseudo quasi-Einstein and also generalized quasi-Einstein, 

        \item it is an Ein$(3)$ spacetime as $S^3+(-\frac{q^2}{r^4})S^2+(-\frac{q^4}{r^8})S+(\frac{q^6}{r^{12}})g=0$,

        \item its nature is of generalized Roter type,

        \item its Ricci tensor is compatible with Riemann, conformal, and projective curvature tensor.

    \end{enumerate}

\end{cor}
If $\lambda=0=q(t)$ and the mass of the body is time-independent in (\ref{VBdS}), then VBdS spacetime turns into Vaidya spacetime. Shaikh et al. \cite{SKS19} investigated the curvature properties of Vaidya spacetime.

\begin{cor}
	The Vaidya spacetime \cite{SKS19} obeys the following curvature  properties:
    \begin{enumerate}[label=(\roman*)]
        \item its scalar curvature vanishes and hence $cir(R)=R$ and $har(R)=C$,

        \item it admits $C\cdot C=\frac{m}{r^3}Q(g,C)$,

        \item satisfies the curvature conditions $R\cdot R-\frac{m}{r^3} Q(g,C)=Q(S,R)$ and \\ $R\cdot C+C\cdot R=Q(S,C)+\frac{2m}{r^3}Q(g,C)$,

        \item its conformal $2$-forms are recurrent for  $\Gamma=\left\lbrace \frac{m}{m'}, 0,0,0\right\rbrace $,

        \item it is Ricci simple such that $S=\varphi (\eta \otimes \eta)$ holds for $\varphi=2m'$ and $\eta=\left\lbrace \frac{1}{r}, 0,0,0\right\rbrace $,

        \item its Ricci tensor is compatible with Rieman and conformal.
    \end{enumerate}
\end{cor}

Again, if $m(t)=$ constant and $\lambda=0=q(t)$, the VBdS spacetime (\ref{VBdS}) reduces to the Schwarzschild spacetime. Thus we have the following:

\begin{cor}
	The curvature properties of the Schwarzschild spacetime are given by the following:
    \begin{enumerate}[label=(\roman*)]
        \item it is Ricci flat and hence $R=cir(R)=C=har(R)=P$,
        
        \item it fulfills $R\cdot R=\frac{m}{r^3} Q(g, R)$ and hence it is a Deszcz pseudosymmetric manifold,
        
        \item it satisfies $divR=0$ and hence the  curvature of the spacetime is harmonic,
        
        \item the compatible tensor for $R$ possesses the following general form: 
	   $$	\left(
	   \begin{array}{cccc}
	   \mathcal{H}_{11} & \mathcal{H}_{12} & 0 & 0 \\
	   \mathcal{H}_{12} & \mathcal{H}_{22} & 0 & 0 \\
	   0 & 0 & \mathcal{H}_{33} & \mathcal{H}_{34} \\
	   0 & 0 & \mathcal{H}_{34} & \mathcal{H}_{44}
	\end{array}
	\right)
	$$
	where $\mathcal{H}_{ij}$ being arbitrary scalars.
    \end{enumerate}
\end{cor}

\section{\bf Symmetry and Ricci soliton admitted on VBdS spacetime}

The set $\mathcal{K}(\mathscr{V}_n)$ of all  Killing vector fields on $\mathscr{V}_n$ constitutes a Lie subalgebra of the Lie algebra of all smooth vector fields. Furthermore, if  $\mathcal{K}(\mathscr{V}_n)$ consists of at most $n(n+1)\slash 2$ linearly independent Killing vector fields, then $\mathscr{V}_n$ is known as a maximally symmetric space i.e., if the cardinality of a basis of $\mathcal{K}(\mathscr{V}_n)$ is exactly the possibly largest number $n(n+1)\slash 2$. If $\mathscr{V}_n$ is maximally symmetric, then the scalar curvature is constant. But the converse is not true because one can take the Schwarzschild spacetime as a counter-example. We mention that in Section 3, we have obtained the scalar curvature of the VBdS spacetime as $\bar\kappa=4\lambda$. The vector field $\frac{\partial}{\partial \phi}$  is Killing and the vector fields $\frac{\partial}{\partial t}$, $\frac{\partial}{\partial r}$, $\frac{\partial}{\partial \theta}$ are non-Killing on the VBdS spacetime. The non-zero components of $\pounds_\frac{\partial}{\partial t}g$ and $\pounds_\frac{\partial}{\partial r}g$ are given by

\begin{eqnarray}
    &&\left(\pounds_{\frac{\partial}{\partial t}}g\right)_{11}=\frac{-2rm'+(q^2)'}{r^2},\ \ \left(\pounds_{\frac{\partial}{\partial t}}g\right)_{22}=0,\ \ \left(\pounds_{\frac{\partial}{\partial t}}g\right)_{33}=0,\ \ \left(\pounds_{\frac{\partial}{\partial t}}g\right)_{44}=0,\notag\\
    &&\left(\pounds_{\frac{\partial}{\partial r}}g\right)_{11}=\frac{-2(r^4 \lambda-3rm+3q^2)}{3r^3},\ \ \left(\pounds_{\frac{\partial}{\partial r}}g\right)_{22}=0,\ \ \left(\pounds_{\frac{\partial}{\partial r}}g\right)_{33}=-2r,\ \ \left(\pounds_{\frac{\partial}{\partial r}}g\right)_{44}=-2r \sin^2\theta.\notag
\end{eqnarray}

In view of the above calculation, we can declare the following:
\begin{thm}
   If $(q^2)'-2rm'> 0,$ then the VBdS spacetime possesses almost $\eta$-Yamabe soliton for the soliton vector field $\xi=\frac{\partial}{\partial t}$ as it satisfies the relation
    $$\frac{1}{2}\pounds_\xi g + (\mu-\bar\kappa)S+ \frac{1}{2}\eta\wedge \eta=0,$$ 
    for $\mu=4\lambda$ and the $1$-form $\eta=\left(-\frac{1}{r}\sqrt{(q^2)'-2rm'},0,0,0\right)$. 
\end{thm}

\begin{thm}
    If $6q^2-2r^7-6rmq^2-6r^4m'+3r^3(q^2)'=0$, then the VBdS spacetime admits almost Ricci soliton for the soliton vector vector field $\xi=\frac{\partial}{\partial r}$ given as follows:    
    $$\frac{1}{2}\pounds_\xi g-\frac{r^3}{2q^2}S-\frac{r^7+q^4}{2rq^4}g=0.$$
\end{thm}

If $\mathcal N=\pounds_\xi har(R)$, then the non-zero components of $\pounds_\xi har(R)$ are given by
\begin{eqnarray}
    \mathcal{N}_{1414}&=&\frac{(2r^4\lambda+3rm-3q^2)(\lambda_1-3r^2)\sin2\theta}{9r^4}=-\mathcal{N}_{1441}=-\mathcal{N}_{4114}=\mathcal{N}_{4141},\notag\\
    \mathcal{N}_{1424}&=&\frac{(2r^4\lambda+3rm-3q^2)\sin2\theta}{3r^2}=-\mathcal{N}_{1442}=\mathcal{N}_{2414}=-\mathcal{N}_{2441}=-\mathcal{N}_{4124}==-\mathcal{N}_{4214},\notag\\
    \mathcal{N}_{3434}&=&\frac{2}{3}\lambda_2 \sin2\theta=-\mathcal{N}_{3443}=-\mathcal{N}_{4334}=\mathcal{N}_{4343} \notag;  
\end{eqnarray}
This leads to the following:
\begin{thm}
    The VBdS spacetime with $q^2-rm\neq0$ admits generalized conharmonic curvature inheritance for the non-Killing vector field $\xi=\frac{\partial}{\partial \theta}$ as it satisfies the relation
    $$\pounds_\xi har(R)=\zeta_1 har(R) + \zeta_2 (g\wedge g) + \zeta_3 (g\wedge S) + \zeta_4 (S\wedge S),$$
where $\zeta_i$ (i=1,2,3,4) are given by
\begin{eqnarray}
    \zeta_1 &=& \frac{2\cos \theta}{\sin ^2\theta},\nonumber\\
    \zeta_2 &=& -\frac{3 \cos\theta(rm-q^2)(3rm-5q^2)^2}{16r^4q^4\sin^2\theta},\nonumber\\
    \zeta_3 &=& -\frac{3\cos\theta(rm-q^2)(3rm-5q^2)}{4q^4\sin^2\theta},\nonumber\\
    \zeta_4 &=& -\frac{3r^4\cos\theta (rm-q^2)}{4q^4\sin^2\theta}. \nonumber   
\end{eqnarray}
\end{thm}

\begin{cor}
    If  $rm=q^2$, then the VBdS spacetime possesses conharmonic curvature inheritance for the vector field $\xi=\frac{\partial}{\partial \theta}$ as it realizes the condition
    $$\pounds_\xi har(R)=\frac{2\cos \theta}{\sin ^2\theta}  har(R) .$$
\end{cor}

\section{\bf Energy momentum tensor of VBdS spacetime}

In general theory of relativity, the energy momentum tensor $T$ of a spacetime is expressed as $$S-\frac{\bar{\kappa}}{2}g+\Lambda g= \frac{8\pi G}{c^4}T$$ 
where $T$ represents the energy momentum tensor, $\bar{\kappa}$ and $S$ are respectively scalar curvature, Ricci curvature of the spacetime. The cosmological constant is denoted by $\Lambda$, the gravitational constant by $G$ and $c$ is the speed of light in a vacuum.

\indent Assuming that $\frac{8\pi G}{c^4}=1$, the energy momentum tensor are computed as given below: 
$$\begin{array}{c}
 T_{11}=\frac{(r^6\lambda(-3+r^2\lambda)-r^2(3+2r^2\lambda)q^2-3q^4+6rm(r^4\lambda+q^2)+3r^3(-2rm'+(q^2)'))}{3r^6}, \\
 T_{12}=\lambda+\frac{q^2}{r^4}, \ \
 T_{33}=\frac{r^4\lambda-q^2}{r^2}, \\
 T_{44}=\frac{(r^4\lambda-q^2)\sin^2\theta}{r^2}.
\end{array}$$ 

 Also, the non-zero components of the tensor $Q(T,R)$ are computed as given below:
  
$$\begin{array}{c}
Q(T,R)_{1313,12}=\frac{(5r^4\lambda-6rm+9q^2)(2rm'-(q^2)')}{3r^5}=\frac{1}{\sin^2\theta}Q(T,R)_{1414,12}, \\

Q(T,R)_{1213,12}=-\frac{(5r^4\lambda-6rm+9q^2)(2rm'-(q^2)')}{6r^5}=\frac{1}{\sin^2\theta}Q(T,R)_{1214,14}, \\

Q(T,R)_{1223,13}=\frac{(9r^5\lambda m-r(14r^3\lambda+3m)q^2+6q^4)}{3r^6}=\frac{1}{\sin^2\theta}Q(T,R)_{1224,14}, \\

Q(T,R)_{2434,13}=\frac{(4r^3\lambda q^2-3m(3r^4\lambda+q^2))\sin^2\theta}{8r^3}=-Q(T,R)_{2334,14}. \\
\end{array}$$ 
From the above, we may conclude the following:

\begin{thm}
The energy momentum tensor of the VBdS spacetime \eqref{VBdS} realizes the following properties:
\begin{enumerate}[label=(\roman*)]
	\item $Q(T,R)=-2\lambda Q(g,R)+Q(S,R),$

	\item the energy momentum tensor $T$ is compatible for the Riemann $(R)$, projective $(P)$, conformal $(C)$, concircular $(cir)$ and conharmonic $(har)$ curvature tensors.	 
\end{enumerate}
\end{thm}
\section{\bf VBdS spacetime Vs. Vaidya-Bonner spacetime}

VBdS spacetime extends the Vaidya-Bonner spacetime which includes radiation as an exterior field. The following is a comparison of the VBdS spacetime with the Vaidya-Bonner spacetime in terms of several types of geometric structures, various kinds of symmetry and the Ricci solitons: 

\begin{enumerate}[label=(\roman*)]
	\item The scalar curvature of Vaidya-Bonner spacetime vanishes while it does not vanish for the VBdS spacetime.

    \item The Ricci tensor is compatible with $R$, $C$ of Vaidya-Bonner spacetime,  while the Ricci tensor is compatible with $R$, $C$, $har(R)$, $P$, $cir(R)$ of VBdS spacetime.
    \item The VBdS spacetime admits almost Ricci soliton for the non-Killing vector field $\frac{\partial}{\partial r}$ while it does not admit by Vaidya-Bonner spacetime. 
    
   \item The VBdS spacetime does not admit almost $\eta$-Yamabe soliton for the non-Killing vector field $\frac{\partial}{\partial \theta}$ but the Vaidya-Bonner spacetime possesses the relation for the same vector field.  
      \item  The VBdS spacetime admits generalized conharmonic curvature inheritance with respect to the non-Killing vector field $\frac{\partial}{\partial \theta}$ whereas the Vaidya-Bonner spacetime does not reveal such inheritance,
      \item For the VBdS spacetime the energy momentum tensor fulfills $Q(T,R)=-2\lambda Q(g,R)+Q(S,R),$ but in the Vaidya-Bonner spacetime the energy momentum tensor satisfies $Q(T,R)=Q(S,R)$. 
     	 
\end{enumerate}
However, they have the following similar properties:
\begin{enumerate}[label=(\roman*)]
    \item both spacetimes are pseudosymmetric due to conformal and conharmonic curvatures,

    \item both of them are generalized Roter type, 

    \item both spacetimes describes Ein$(3)$ and $2$-quasi-Einstein manifolds,

    \item conformal curvature $2$-forms are recurrent for both spacetimes,
     \item the vector field $\frac{\partial}{\partial \phi}$ is Killing in both spacetimes,
        \item the vector fields $\frac{\partial}{\partial t}$, $\frac{\partial}{\partial r}$ and  $\frac{\partial}{\partial \theta}$ are non-Killing vector fields in both spacetimes, 
        \item both the spacetimes possesses almost $\eta$-Yamabe soliton for the non-Killing  vector field $\frac{\partial}{\partial t}$,
        \item the energy momentum tensor is compatible for Riemann, projective, conformal, concircular and conharmonic  curvatures for both the spacetimes.
\end{enumerate}

\vspace{1em}
\section{\bf  Acknowledgment}
 The third author greatly acknowledges to The University Grants Commission, Government of India for the award of Senior Research Fellow. All the algebraic computations of Section $3$ to $4$ are performed by a program in Wolfram Mathematica developed by the first author (A. A. Shaikh).

	{\bf Data Availability:}- No data was used for the research described in the article.

{\bf Conflict of Interest:}-There is no conflict of interests regarding publication of this paper.


%



\begin{thebibliography}{99}\baselineskip=16pt


\bibitem{AD83}
Adam\'{o}w, A. and Deszcz, R.,
\emph{On totally umbilical submanifolds of some class of Riemannian manifolds},
Demonstratio Math.,
\textbf{16} (1983), 39--59.


\bibitem{Ahsan1977_231}
Ahsan, Z.,
\emph{Algebraic classification of space-matter tensor in general relativity},
Indian J. Pure Appl. Math.,
\textbf{8(2)} (1977), 231--237


\bibitem{Ahsan1977_1055}
Ahsan, Z.,
\emph{Algebra of space-matter tensor in general relativity},
Indian J. Pure Appl. Math.,
\textbf{8(9)} (1977), 1055--1061.



\bibitem{Ahsan1978}
Ahsan, Z.,
\emph{Collineation in electromagnetic field in general relativity- The null field
	case},
Tamkang J. Maths.,
\textbf{9(2)} (1978), 237.


\bibitem{Ahsan1987} 
Ahsan, Z.,
\emph{On the Nijenhuis tensor for null electromagnetic field},
J. Math. Phys. Sci.,
\textbf{21(5)} (1987), 515--526.


\bibitem{Ahsan1995}
Ahsan, Z.,
\emph{Symmetries of the electromagnetic fields in general relativity},
Acta Phys. Sinica,
\textbf{4} (1995), 337.


\bibitem{Ahsan1996}
Ahsan, Z.,
\emph{A symmetry property of the space-time of general relativity in terms of the
	space-matter tensor},
Braz. J. Phys. 
\textbf{26(3)} (1996), 572-576.



\bibitem{Ahasan2005}
Ahsan, Z.,
\emph{On a geometrical symmetry of the space-time of General Relativity},
Bull. Call. Math. Soc.,
\textbf{97(3)} (2005), 191.

\bibitem{Ahsan2018}
Ahsan, Z.,
\emph{Ricci solitons and the spacetime of general relativity},
J. Tensor Soc.,
\textbf{12} (2018), 49--64.

\bibitem{AA2012}
Ahsan, Z. and Ali, M.,
\emph{Symmetries of type D pure radiation fields}.
Int. J. Theo. Phys.,
\textbf{51} (2012), 2044-2055.

\bibitem{AhsanAli2014}
Ahsan, Z. and Ali, M.,
\emph{On some properties of $W$-curvature tensor},
Palestine J. Math.,
\textbf{3(1)} (2014), 61--69.
\bibitem{AH1980}
Ahsan, Z. and Husain, S. I.,
\emph{Null electromagnetic fields, total gravitational radiation
	and collineations in general relativity},
Annali di Mathematical Pura ed Applicata,
\textbf{126} (1980), 379396.
\bibitem{AliAhsan2012}
Ali, M. and Ahsan Z.,
\emph{Ricci solitons and symmetries of spacetime manifold of general relativity}
Global J. Adv. Research Classical Mod. Geom., 
\textbf{1(2)} (2012), 75--84.
\bibitem{AliAhsan2013}
Ali, M. and  Ahsan, Z.,
\emph{Geometry of Schwarzschild soliton},
J. Tensor Soc.,
\textbf{7} (2013), 49--57.
\bibitem{AliAhsan2015}
Ali, M. and Ahsan, Z.,
\emph{Gravitational field of Schwarzschild soliton},
Arab J. Math. Sci.,
\textbf{21(1)} (2015), 15--21.





\bibitem{ARS95}
Al\'ias, L. J., Romero, A. and S\'anchez, M., \emph{Uniqueness of complete spacelike hypersurfaces of constant mean curvature in generalized Robertson-Walker space-times}, Gen. Relativity Gravitation, \textbf{27(1)} (1995), 71--84.



\bibitem{ADEHM14}
Arslan, K., Deszcz, R., Ezenta\c{s}, R., Hotlo\'{s}, M. and Murathan, C., \emph{On generalized Robertson-Walker spacetimes satisfying some curvature condition}, Turkish J. Math., \textbf{38(2)} (2014), 353--373.
\bibitem{Bess87}
Besse, A. L.,
\emph{Einstein Manifolds},
Springer-Verlag, Berlin, Heidelberg, \textbf{1987}.

\bibitem{Blaga2016}
Blaga, A. M.,
\emph{$eta$-Ricci solitons on Lorentzian para-Sasakian manifolds},
Filomat, 
\textbf{30(2)} (2016), 489--496.
\bibitem{Brink1925}
Brinkmann, H. W., \emph{Einstein spaces which are mapped conformally on each other}, Math. Ann. \textbf{94}
(1925), 119--145.













\bibitem{Cart26}
Cartan, \'E., \emph{Sur une classe remarquable d'espaces de Riemannian}, Bull. Soc. Math. France, \textbf{54} (1926), 214- 264.

\bibitem{Cart46}
Cartan, \'E.,
\emph{Le\c cons sur la g\' eom\' etrie des espaces de Riemann}, 2nd ed., Paris, \textbf{1946}.

\bibitem{Chak87}
Chaki, M. C.,
\emph{On pseudosymmetric manifolds},
An. \c{S}tiin\c{t}.  Univ. AL. I. Cuza Ia\c{s}i. Mat. (N.S.)  Sect. Ia,
\textbf{33(1)} (1987), 53--58.

\bibitem{Chak88}
Chaki, M. C., \emph{On pseudo Ricci symmetric manifolds}, Bulgarian J. Phys., {\bf 15} (1988), 526--531.

\bibitem{VBdS_2007_Hawking}
Chen, D. and Yang, S.,
\emph{Hawking radiation of the Vaidya–Bonner–de Sitter black hole}, 
New J. Phys, \textbf{9(8)} (2007), 252.


\bibitem{Cho2009}
Cho, J. and Kimura, M.,
\emph{Ricci solitons and real hypersurfaces in a complex space form}, 
Tohoku Math. J.,
\textbf{61(2)} (2009), 205--212.







\bibitem{DDHKS00}
Defever, F.,  Deszcz, R.,  Hotlo$\acute{\mbox{s}}$, M.,  Kucharski, M.  and  Sent$\ddot{\mbox{u}}$rk, Z.,  \emph{Generalisations of Robertson-Walker spaces}, Ann. Univ. Sci. Budapest, E$\ddot{\mbox{o}}$tv$\ddot{\mbox{o}}$s Sect. Math., {\bf 43} (2000), 13--24.
\bibitem{DD91}
Defever, F. and  Deszcz, R., \emph{On semi-Riemannian manifolds satisfying the condition $R\cdot R=Q(S,R)$}, in:Geometry and Topology of Submanifolds III, World Sci., River Edge, NJ, (1991), 108-130.


\bibitem{Desz92}
Deszcz, R., \emph{On pseudosymmetric spaces},
Bull. Belg. Math. Soc., Ser. A,
\textbf{44} (1992), 1--34.

\bibitem{Desz93}
Deszcz, R., \emph{Curvature properties of a pseudosymmetric manifolds}, Colloq. Math., \textbf{62} (1993), 139--147.

\bibitem{Desz03}
Deszcz, R., \emph{On Roter type manifolds}, 5-th Conference on Geometry and Topology of Manifolds, Krynica, Poland, April 27 - May 3, (2003),  25.

\bibitem{DG02}
Deszcz, R. and G\l ogowska, M.,
\emph{Some examples of nonsemisymmetric Ricci-semisymmetric hypersurfaces},
Colloq. Math., \textbf{94} (2002), 87--101.

\bibitem{DGHS98}
Deszcz, R., G\l ogowska, M., Hotlo\'{s}, M. and \d Sent\"{u}rk, Z.,
\emph{On certain quasi-Einstein semi-symmetric hypersurfaces},
Ann. Univ. Sci. Budapest E\"{o}tv\"{o}s Sect. Math.,
\textbf{41} (1998), 151--164.



\bibitem{DGHS11}
Deszcz, R., G\l ogowska, M., Hotlo\'s, M. and Sawicz, K.,
\emph{A survey on generalized Einstein metric conditions},
Advances in Lorentzian Geometry, Proceedings of the Lorentzian Geometry Conference in Berlin, AMS/IP Studies in Advanced Mathematics, \textbf{49}, S.-T. Yau (series ed.),
M. Plaue, A.D. Rendall and M. Scherfner (eds.), 2011, 27-46.

\bibitem{DGHZ15} 
Deszcz, R., G\l ogowska, M., Hotlo\'s, M. and Zafindratafa, G.,
\emph{On some curvature conditions of pseudosymmetric type}, 
Period. Math. Hungarica,
\textbf{70(2)} (2015), 153--170.

\bibitem{DGHZ16}
Deszcz, R., G\l ogowska, M., Hotlo\'s, M. and Zafindratafa, G., \emph{Hypersurfaces in space
forms satisfying some curvature conditions}, J. Geom. Phys., \textbf{99} (2016), 218--231.
\bibitem{DGJPZ13}
Deszcz, R., G\l ogowska, M., Je\l owicki, L., Petrovi\'{c}-Torga\u{s}ev, M. and Zafindratafa, G.,
\emph{On Riemann and Weyl compatible tensors},
Publ. Inst. Math. (Beograd) (N.S.),
\textbf{94(108)} (2013), 111--124.
\bibitem{DGJZ-2016}
Deszcz, R., G\l ogowska, M., Je\l owicki, J. and Zafindratafa, Z.,
\emph{Curvature properties of some class of warped product manifolds},
Int. J. Geom. Methods Mod. Phys.,
\textbf{13} (2016), 1550135.

\bibitem{DGP-TV-2015}
Deszcz, R., G\l ogowska, M., Petrovi\'{c}-Torga\u{s}ev, M. and Verstraelen, L.,
\emph{Curvature properties of some class of minimal hypersurfaces in Euclidean spaces},
Filomat,
\textbf{29} (2015), 479--492.

\bibitem{DGPSS11}
Deszcz, R., G\l ogowska, M., Plaue, M., Sawicz, K. and Scherfner, M.,
\emph{On hypersurfaces in space forms satisfying particular curvature conditions of Tachibana type},
Kragujevac J. Math.,
\textbf{35} (2011), 223--247.

\bibitem{DGP-TV-2011}
Deszcz, R., G\l ogowska, M., Petrovi\'c-Torga\u{s}ev, M. and Verstraelen, L.,
\emph{On the Roter type of Chen ideal submanifolds},
Results Math.,
\textbf{59} (2011), 401--413.


\bibitem{DH03}
Deszcz, R. and Hotlo\'{s}, M.,
\emph{On hypersurfaces with type number two in spaces of constant curvature},
Ann. Univ. Sci. Budapest E\"{o}tv\"{o}s Sect. Math.,
\textbf{46} (2003), 19--34.

\bibitem{DHJKS14}
Deszcz, R., Hotlo\'{s}, M., Je\l owicki, J., Kundu, H. and Shaikh, A. A.,
\emph{Curvature properties of G\"{o}del metric},
Int. J. Geom. Methods Mod. Phys., \textbf{11} (2014), 1450025. Erratum: \emph{Curvature properties of Gödel metric}, Int. J. Geom. Methods Mod. Phys., \textbf{16} (2019), 1992002.

\bibitem{DK99}
Deszcz, R. and Kucharski, M., \emph{On curvature properties of certain generalized Robertson-Walker spacetimes}, Tsukuba J. Math., \textbf{23(1)} (1999), 113--130.

\bibitem{DPSch-2013}
Deszcz, R., Plaue, M. and Scherfner, M.,
\emph{On Roter type warped products with 1-dimensional fibres},
J. Geom. Phys., \textbf{69} (2013), 1--11.
\bibitem{Duggal1992}
Duggal, K. L., 
\emph{Curvature inheritance symmetry in Riemannian spaces with applications to fluid space times},
J. Math. Phys., \textbf{33(9)} (1992), 2989--2997.

\bibitem{EC21}
Eyasmin, S. and Chakraborty, D., \emph{Curvature properties of (t-z)-type plane wave metric}, J. Geom. Phys., \textbf{160} (2021), 104004.
\bibitem{ECS22}
Eyasmin, S., Chakraborty, D. and Sarkar, M., \emph{Curvature properties of Morris-Thorne Wormhole metric}, J. Geom. Phys., \textbf{174(2)} (2022), 104457.





\bibitem{EDS_sultana_2022}
Eyasmin, S., Datta, B. R. and Sarkar, M., \emph{On sultana-dyer spacetime: curvatures and geometric structures}. Int. J. Geom. Methods Mod. Phys.  (2022) DOI: 10.1142/S0219887823501013


\bibitem{VBdS_2012}
Farmany, A.,  
\emph{Thermodynamics of Vaidya-Bonner-de Sitter space time based on the generalized space-time uncertainty},
Astrophys. Space Sci, \textbf{337} (2012), 785--787.
\bibitem{F81}
Ferus, D., \emph{A remark on Codazzi tensors on constant curvature space}, Glob. Diff. Geom. Glob. Ann., Lecture
notes 838, Springer, \textbf{1981}.


\bibitem{Gray78}
Gray, A., \emph{Einstein-like manifolds which are not Einstein}, Geom. Dedicta, \textbf{7} (1978), 259--280.

\bibitem{GP09}
Griffiths, J. B. and Podolsk\'y, J., \emph{Exact space-times in Einstein’s general relativity}, Cambridge University Press, \textbf{2009}.




\bibitem{Glog02}
G\l ogowska, M.,
\emph{Semi-Riemannian manifolds whose Weyl tensor is a Kulkarni-Nomizu square},
Publ. Inst. Math. (Beograd) (N.S.),
\textbf{72(86)} (2002), 95--106.

\bibitem{Glog-2007}
G\l ogowska, M., \emph{On Roter type manifolds}, Pure and Applied Differential Geometry- PADGE, (2007), 114--122.
	\bibitem{G08}
	G\l ogowska, M., \emph{On quasi-Einstein Cartan type hypersurfaces}, J. Geom. Phys. \textbf{58} (2008), 599--614.

\bibitem{Guler2019}
G\"uler, S. and Cr\'{a}\c{s}mare\v{a}nu, M., \emph{Ricci-Yamabe maps for Riemannian flow and their volume
variation and volume entropy}, Turk. J. Math., \textbf{43}  (2019), 2631--2641.
\bibitem{HV07}
Haesen, S. and Verstraelen, L., \emph{Properties of a scalar curvature invariant depending on two planes}, Manuscripta Math., \textbf{122} (2007), 59-72.
\bibitem{HV09}
Haesen, S. and Verstraelen, L., \emph{Natural intrinsic geometrical symmetries}, Symmetry, Integrability and Geometry, Methods and Appl. SIGMA, \textbf{5} (2009), 086, 15 pages.
\bibitem{Hamilton1982}
Hamilton, R. S.,
\emph{Three manifolds with positive Ricci curvature},
J. Diff. Geom., 
\textbf{17} (1982), 255--306.


\bibitem{Hamilton1988}
Hamilton, R. S.,
\emph{The Ricci flow on surfaces},
Contemp. Math.,
\textbf{71} (1988), 237--261.













\bibitem{KLD1969}
Katzin, G. H., Livine, J. and Davis, W. R., 
\emph{Curvature collineations: A fundamental symmetry property of the space-times of general relativity defined by the vanishing Lie derivative of the Riemann curvature tensor},
J. Math. Phys., 
\textbf{10(4)} (1969), 617--629.



\bibitem{KLD1970}
Katzin, G. H., Livine, J. and Davis, W. R.,
\emph{Groups of curvature collineations in Riemannian space-times which admit fields of parallel vectors},
J. Math. Phys.,
\textbf{11} (1970), 1578--1580.


%

\bibitem{Kowa06}
Kowalczyk, D., \emph{On the Reissner-Nordström-de Sitter type spacetimes}, Tsukuba J. Math., \textbf{30(2)} (2006), 363--381.

\bibitem{VBdS_1999}
Li, Zhong-heng; You, Liang; Mi, Li-qin,
\emph{New quantum effect for Vaidya-Bonner-de Sitter black holes},
Int. J. Theor. Phys. \textbf{38(3)} (1999), 925--931.



\bibitem{LR89}
Lovelock, D. and Rund, H., \emph{Tensors, differential forms and variational principles}, Courier Dover Publications, \textbf{1989}.

\bibitem{MM12a}
Mantica, C. A. and Molinari, L. G., \emph{Extended Derdzinski-Shen theorem for curvature tensors}, Colloq. Math., \textbf{128} (2012), 1--6.

\bibitem{MM22a}
Mantica, C. A. and Molinari, L. G., \emph{Spherical doubly warped spacetimes for radiating stars and cosmology}, Gen. Relativ. Gravit., \textbf{54}: 98 (2022) (34 pages). DOI: 10.1007/s10714-022-02984-7
\bibitem{MM22b}
Mantica, C. A. and Molinari, L. G., \emph{The Jordan algebras of Riemann, Weyl and curvature compatible tensors}, Colloq. Math., \textbf{167} (2022), 63--72.


\bibitem{MM12b}
Mantica, C. A. and Molinari, L. G., \emph{Riemann compatible tensors}, Colloq. Math., \textbf{128} (2012), 197--210.
\bibitem{MM13}
Mantica, C. A. and Molinari, L. G., \emph{Weyl compatible tensors}, Int. J. Geom. Methods Mod. Phys., \textbf{11(08)} (2014), 1450070.

\bibitem{MS12a}
Mantica, C. A. and Suh, Y. J., \emph{The closedness of some generalized curvature 2-forms on a Riemannian manifold I}, Publ. Math. Debrecen, {\bf{81(3-4)}} (2012), 313--326.

\bibitem{MS13a}
Mantica, C. A. and Suh, Y. J., \emph{The closedness of some generalized curvature 2-forms on a Riemannian manifold II}, Publ. Math. Debrecen, {\bf{82(1)}} (2013), 163--182.
 
\bibitem{MS14}
Mantica, C. A. and Suh, Y. J., \emph{Recurrent conformal 2-forms on pseudo-Riemannian manifolds}, Int. J. Geom. Methods Mod. Phy., \textbf{11(6)} (2014), 1450056 (29 pages).





\bibitem{MIKES76} 
Mike$\check{s}$, J., \emph{Geodesic mappings of symmetric Riemannian spaces}, Odessk. Univ., \textbf{3924-76} (1976), 1-10.




\bibitem{MIKES88}
Mike\v s, J., \emph{Geodesic mappings of special Riemannian spaces}, Topics in differential
geometry, Vol. I, II (Debrecen, 1984), Colloq. Math. Soc. J\'anos Bolyai, 46, North-Holland, Amsterdam, (1988), 793--813.


\bibitem{MIKES92}
Mike\v s, J.,\emph{Geodesic mappings of m-symmetric and generalized semisymmetric spaces}
(Russian), translated from Izv. Vyssh. Uchebn. Zaved. Mat. 1992, no. 8, 42-46 Russian Math. (Iz. VUZ), \textbf{36} (1992), no. 8, 38--42.

\bibitem{MIKES96} 
Mike$\check{s}$, J., \emph{Geodesic mapping of affine-connected and Riemannian spaces}, J. Math. Sci., \textbf{78(3)} (1996), 311--333. 
\bibitem{MS94}
Mike\v s, J. and Sobchuk, V.S., \emph{Geodesic mappings of 3-symmetric Riemannian spaces}
(Russian), translated from Ukrain. Geom. Sb. no. 34 (1991), 80--83, iii J. Math. Sci.
69 (1994), no. 1, 885--887.
\bibitem{MSV15}
Mike\v s, J., Stepanova, E. and Van\v zurov\'a, A., \emph{ Differential geometry of special mappings},
Palack\'y University Olomouc, Faculty of Science, Olomouc (2015), 568 pp.

\bibitem{MVH09}
Mike\v{s}, J., Van\v{z}urov\'{a}, A. and Hinterleitner, I, \emph{Geodesic mappings and some generalizations}, Palacky Univ. Press, Olomouc, \textbf{2009}.









\bibitem{VBdS_2013}
Pan, W. Z., Yang, X. J. and Yu, G. X.,
\emph{Quantum nonthermal effect of the Vaidya–Bonner–de Sitter black hole},
Front. Phys. \textbf{9(1)} (2014), 94–97.
\bibitem{VBdS_1987}
Patino, A. and Rago, H., \emph{A radiating charge embedded in a De Sitter universe},   Phys. Lett. A, \textbf{121(7)} (1987), 329-330.

\bibitem{Patt52} 
Patterson, E. M., \emph{Some theorems on Ricci recurrent spaces}, J. London Math. Soc., \textbf{27} (1952), 287--295.



\bibitem{Pigola2011}
Pigola, S., Rigoli, M., Rimoldi, M., Setti, A. G., 
\emph{Ricci almost solitons},
Ann. Scuola Norm. Sup. Pisa Cl. Sci., \textbf{X(5)} (2011), 757--799.




\bibitem{P95}
Prvanovi$\acute{\mbox{c}}$,  M.,  \emph{On weakly symmetric Riemannian manifolds,} Publ. Math. Debrecen, \textbf{46(1-2)} (1995), 19--25.


\bibitem{Ruse46}
Ruse, H. S., \emph{On simply harmonic spaces}, J. London Math. Soc., \textbf{21} (1946), 243--247.
\bibitem{Ruse49a}
Ruse, H. S., \emph{On simply harmonic `kappa spaces' of four dimensions}, Proc. London Math. Soc., \textbf{50}
(1949), 317--329.
\bibitem{Ruse49b}
Ruse, H. S., \emph{Three dimensional spaces of recurrent curvature}, Proc. London Math. Soc., \textbf{50} (1949),
438--446.




\bibitem{S09}
Shaikh, A. A., \emph{On pseudo quasi-Einstein manifolds}, Period. Math. Hungarica, \textbf{59(2)} (2009), 119--146.

\bibitem{SAA18}
Shaikh, A. A., Ali, M. and Ahsan, Z., \emph{Curvature properties of Robinson-Trautman metric}, J. Geom., \textbf{109(2)} (2018), 1--20. DOI: 10.1007/s00022-018-0443-1



\bibitem{SAAC20}
Shaikh, A. A., Ali, A., Alkhaldi, A. H. and Chakraborty, D., \emph{Curvature properties of Melvin magnetic metric}, J. Geom. Phys., \textbf{150} (2020), 103593. DOI: 10.1016/j.geomphys.2019.103593.

\bibitem{SAAC20N}
Shaikh, A. A., Ali, A., Alkhaldi, A. H. and Chakraborty, D., \emph{Curvature properties of Nariai spacetimes}, Int. J. Geom. Methods Mod. Phys., \textbf{17(03)} (2020), 2050034. DOI: 10.1142/S0219887820500346

\bibitem{SAACD_LTB_2022}
Shaikh, A. A., Ali, A., Alkhaldi, A. H., Chakraborty, D. and Datta, B. R., \emph{On some curvature properties of Lemaitre–Tolman–Bondi spacetime}. Gen. Relativ. Gravit. \textbf{54(1)} (2022), 6 (21 pages). DOI: 10.1007/s10714-021-02890-4
\bibitem{SAD_pgm_2023}
Shaikh, A. A., Ahmed, F. and Datta, B. R.,
\emph{Geometrical properties of a point-like global monopole spacetime}, arXiv:2301.04897 (2023).	

\bibitem{SAR13}
Shaikh, A. A., Al-Solamy, F. R. and Roy, I., 
\emph{On the existence of a new class of semi-Riemannian manifolds}, 
Mathematical Sciences, \textbf{7} (2013), 46.



\bibitem{SB08}
Shaikh, A. A. and Binh, T. Q., \emph{On some class of Riemannian manifolds,} Bull. Transilv. Univ., \textbf{15(50)} (2008), 351--362.
	\bibitem{SBK21}
		Shaikh, A. A., Binh, T. Q. and Kundu, H., \emph{Curvature properties of generalized pp-wave metric}, Kragujevac J. Math., \textbf{45(2)} (2021), 237--258.

\bibitem{SC21}
Shaikh, A. A. and Chakraborty, D., \emph{Curvature properties of Kantowski-Sachs metric}, J. Geom. Phys., \textbf{160} (2021), 103970. DOI: 10.1016/j.geomphys.2020.103970


\bibitem{SDAA_LCS_2021}
Shaikh, A. A., Datta, B. R., Ali, A. and Alkhaldi, A. H., \emph{LCS-manifolds and Ricci solitons}. Int. J. Geom. Methods Mod. Phys., \textbf{18(09)} (2021), 2150138.

\bibitem{ShaikhDatta2022} Shaikh, A. A. and Datta, B. R., \emph{Ricci solitons and curvature inheritance on Robinson-Trautman spacetimes} (2022) arXiv preprint arXiv:2209.03749.



\bibitem{SDC}
Shaikh, A. A., Datta, B. R. and Chakraborty, D., \emph{On some curvature properties of Vaidya-Bonner metric}, Int. J. Geom. Methods. Phys., https://doi.org/10.1142/S0219887821502054.


\bibitem{SDKC19}
Shaikh, A. A., Das, L., Kundu, H. and Chakraborty, D., \emph{Curvature properties of Siklos metric}, Diff.
Goem.- Dyn. Syst., \textbf{21} (2019), 167--180.




\bibitem{SDHJK15}
Shaikh, A. A., Deszcz, R., Hotlo\'{s}, M., Je\l owicki, J. and Kundu, H.,
\emph{On pseudosymmetric manifolds},
Publ. Math. Debrecen,
\textbf{86(3-4)} (2015), 433-456.

	\bibitem{SHDS_hayward}
		Shaikh, A. A., Hui, S. K., Datta, B. R. and Sarkar, M.,
		\emph{On curvature related geometric properties of Hayward black hole spacetime}, (2023), 	arXiv:2303.00932.
		

\bibitem{SHS_warped}
		Shaikh, A. A., Hui, S. K. and Sarkar, M., \emph{Curvature properties of a warped product metric}, Palestine J. Math. (accepted).
		
		\bibitem{SHS_Bardeen}
				Shaikh, A. A., Hui, S. K. and Sarkar, M., \emph{Curvature properties of Bardeen black hole spacetime}, Bulgarian J. Phys., \textbf{50} (2023), 168--189. 



\bibitem{SKH11}
Shaikh, A. A., Kim, Y. H. and Hui, S. K., \emph{On Lorentzian quasi Einstein manifolds}, J. Korean Math. Soc., \textbf{48} (2011), 669--689. Erratum: \emph{On Lorentzian quasi Einstein manifolds,} J. Korean Math. Soc., \textbf{48(6)} (2011), 1327--1328.


\bibitem{SK14}
Shaikh, A. A. and Kundu, H., \emph{On equivalency of various geometric structures}, J. Geom., \textbf{105} (2014), 139--165. DOI: 10.1007/s00022-013-0200-4

\bibitem{SK16}
Shaikh, A. A. and Kundu, H., \emph{On warped product generalized Roter type manifolds}, Balkan J. Geom. Appl., \textbf{21(2)} (2016), 82--95.

\bibitem{SK16srs}
Shaikh, A. A. and Kundu, H., \emph{On curvature properties of Som-Raychaudhuri spacetime}, J. Geom., \textbf{108(2)} (2016), 501--515.
	


\bibitem{SKppsnw}
Shaikh, A. A. and Kundu, H., \emph{On warped product manifolds satisfying some pseudosymmetric type conditions}, Diff. Geom. - Dyn. Syst., \textbf{19} (2017), 119--135.


\bibitem{SK19}
Shaikh, A. A. and Kundu, H., \emph{On generalized Roter type manifolds},
Kragujevac J. Math., \textbf{43(3)} (2019), 471--493.


\bibitem{SKA18}
Shaikh, A. A., Kundu, H. and Ali, Md. S., \emph{On warped product super generalized recurrent manifolds}, An. \c{S}tiin\c{t}. Univ. Al. I. Cuza Ia\c{s}i. Mat. (N. S.), \textbf{LXIV(1)} (2018), 85--99.

\bibitem{SKS19}
Shaikh, A. A., Kundu, H. and Sen, J., \emph{Curvature properties of Vaidya metric}, Indian J. Math., \textbf{61(1)} (2019), 41--59.
\bibitem{VBdS_2022}
Sha, B. and Liu, Z. E., 
\emph{Lorentz-breaking theory and tunneling radiation correction to Vaidya–Bonner de Sitter Black Hole}, 
Eur. Phys. J. C, \textbf{82(7)} (2022), 648.

\bibitem{SP10}
Shaikh, A. A. and Patra, A., \emph{On a generalized class of recurrent manifolds}, 
Arch. Math. (BRNO), {\bf 46} (2010), 71--78.
\bibitem{SR10}
Shaikh, A. A.  and Roy, I.,  \emph{On quasi generalized recurrent manifolds}, Math. Pannon, \textbf{21(2)} (2010), 251--263.
\bibitem{SR11}
Shaikh, A. A. and Roy, I., \emph{On weakly generalized recurrent manifolds}, Ann. Univ. Sci. Budapest, E$\ddot{\mbox{o}}$tv$\ddot{\mbox{o}}$s Sect. Math., \textbf{54} (2011), 35--45.
\bibitem{SRK18}
Shaikh, A. A., Roy, I. and Kundu, H.,
\emph{On the existence of a generalized class of recurrent manifolds},  An. \c{S}tiin\c{t}. Univ. Al. I. Cuza Ia\c{s}i. Mat. (N. S.), \textbf{LXIV(2)} (2018), 233--251.

\bibitem{SRK17}
Shaikh, A. A., Roy, I. and Kundu, H., \emph{On some generalized recurrent manifolds}, Bull.
Iranian Math. Soc., \textbf{43(5)} (2017), 1209--1225.



\bibitem{SSC19}
Shaikh, A. A., Srivastava, S. K. and Chakraborty, D., \emph{Curvature properties of anisotropic scale invariant
metrics}, Int. J. Geom. Meth. Mod. Phys., {\bf 16} (2019), 1950086.

\bibitem{SYH09}
Shaikh, A. A., Yoon, D. W. and Hui, S. K., \emph{ On quasi-Einstein spacetimes}, Tsukuba J. Math., {\bf 33(2)} (2009), 305--326.


\bibitem{Siddiqi2020}
Siddiqi, M. and Akyol, M. A.,  \emph{$\eta $-Ricci-Yamabe soliton on Riemannian submersions from Riemannian manifolds}. arXiv preprint arXiv:2004.14124, (2020).

\bibitem{S81}
Simon, U., \emph{Codazzi tensors}, Glob. Diff. Geom. and Glob. Ann., Lecture notes, 838, Springer-Verlag, 1981,
289--296.

\bibitem{SKP03}
Suh, Y. J., Kwon, J-H. and Pyo, Y. S., \emph{On semi-Riemannian manifolds satisfying the second Bianchi identity}, J. Korean Math. Soc., \textbf{40(1)} (2003), 129--167.
\bibitem{Szab82}
Szab$\acute{\mbox{o}}$, Z. I., \emph{Structure theorems on Riemannian spaces satisfying $R(X,Y)\cdot R=0$, I. The local version}, J. Diff. Geom., \textbf{17} (1982), 531--582.

\bibitem{Szab84}
Szab\'o, Z. I., \emph{Classification and construction of complete hypersurfaces satisfying 
$R(X, Y)\cdot R = 0$}, Acta Sci. Math., \textbf{47} (1984), 321--348.

\bibitem{Szab85}
Szab\'o, Z. I., \emph{Structure theorems 
on Riemannian spaces satisfying $R(X, Y)\cdot R = 0$, II, The global version}, 
Geom. Dedicata, \textbf{19} (1985), 65--108.



\bibitem{Tach74}
Tachibana, S.,
\emph{A theorem on Riemannian manifolds of positive curvature operator},
Proc. Japan Acad.,
\textbf{50} (1974), 301--302.

\bibitem{TB89}
Tam\'{a}ssy, L. and Binh, T. Q., \emph{On weakly symmetric and weakly projective symmetric Riemannian manifolds}, Colloq. Math. Soc. J. Bolyai, \textbf{50} (1989), 663--670.

\bibitem{TB93}
Tam\'{a}ssy, L. and  Binh, T. Q., \emph{On weak symmetries of  Einstein and Sasakian manifolds}, Tensor (N. S.), \textbf{53} (1993), 140--148.




\bibitem{Venz85} 
Venzi, P., 
\emph{Una generalizzazione degli spazi ricorrenti}, 
Rev. Roumaine Math. Pures Appl.,
\textbf{30} (1985), 295--305.

\bibitem{Walk50}
Walker, A. G.,
\emph{On Ruse's spaces of recurrent curvature},
Proc. London Math. Soc.,
\textbf{52} (1950), 36--64.




\bibitem{VBdS_2001}    
Yong, M.,
\emph{Thermal radiation of scalar particles in Vaidya-Bonner-de Sitter space-time},
J. Chongqing Norm. Univ. Nat. Sci. Ed. \textbf{1} (2001), 8--12.



\end{thebibliography}
\end{document}